\documentclass[final]{siamonline1116}

\usepackage{slashbox}
\usepackage[utf8]{inputenc}
\usepackage{bm}
\usepackage{graphicx}
\usepackage{amsmath}
\usepackage[parfill]{parskip}
\usepackage{xcolor}
\usepackage[english]{babel}
\usepackage{amsfonts}
\usepackage{amssymb}
\usepackage{amsbsy}
\usepackage{hyperref}
\usepackage{enumitem}
\usepackage{epstopdf}
\usepackage{subcaption}
\usepackage{cleveref}

\ifpdf
  \DeclareGraphicsExtensions{.eps,.pdf,.png,.jpg}
\else
  \DeclareGraphicsExtensions{.eps}
\fi

\numberwithin{theorem}{section}

\newsiamremark{remark}{Remark}
\newsiamremark{Example}{Example}
\newsiamthm{assumption}{Assumption}

\newcommand{\R}{\mathcal{R}}
\newcommand{\Rs}{\R^d}
\newcommand{\Ro}{\R^s}
\newcommand{\fatQ}{\mathbb{Q}}
\newcommand{\fatR}{\mathbb{R}}
\newcommand{\fatB}{\mathbb{B}}
\newcommand{\fatS}{\mathbb{S}}
\newcommand{\truth}{z}
\newcommand{\est}{x}
\newcommand{\err}{\xi}
\newcommand{\D}[1]{\ensuremath{\mathrm{D}{#1}}}

\DeclareMathOperator{\Rank}{rank}

\newcommand{\matpair}[2]{\left({#1},{#2}\right)}

\newcommand{\TheTitle}{A detectability criterion and data assimilation for non-linear differential equations}
\newcommand{\ShortTitle}{Detectability and data assimilation for non-linear ODEs}
\newcommand{\TheAuthors}{Jason Frank \& Sergiy Zhuk}

\ifpdf
\hypersetup{
  pdftitle={\TheTitle},
  pdfauthor={\TheAuthors}
}
\fi

\headers{\ShortTitle}{\TheAuthors}

\title{{\TheTitle}\thanks{Submitted to the editors \today.
}}

\author{
Jason Frank\thanks{Utrecht University, Mathematical Institute, P.O.~Box 80010, 3508 TA Utrecht, Netherlands,
	(\email{j.e.frank@uu.nl}).}
\and
Sergiy Zhuk\thanks{IBM Research, Dublin, Ireland
	(\email{sergiy.zhuk@ie.ibm.com}).}
}

\usepackage{amsopn}

\begin{document}
\maketitle

\begin{abstract}
In this paper we propose a new sequential data assimilation method for non-linear ordinary differential equations with compact state space. The method is designed so that the Lyapunov exponents of the corresponding estimation error dynamics are negative, i.e.~the estimation error decays exponentially fast. The latter is shown to be the case for generic \emph{regular} flow maps if and only if the observation matrix $H$ satisfies detectability conditions: the rank of $H$ must be at least as great as the number of nonnegative Lyapunov exponents of the underlying attractor. Numerical experiments illustrate the exponential convergence of the method and the sharpness of the theory for the case of Lorenz '96 and Burgers equations with incomplete and noisy observations.
\end{abstract}



\section{Introduction}
Consider a process described by the following ordinary differential equation (ODE):
\begin{equation}\label{model}
	\dot{\truth}  = f(t,\truth), \qquad \truth(t) \in \mathcal{D}\subset \R^d, \quad \truth(0) = \truth_0.
\end{equation}
The state $\truth(t)$ and initial condition $\truth_0$ are presumed to be unknown, but information about the state is obtained via a noisy observation process:
\begin{equation}\label{obs}
	y(t) = H(t) \truth(t) + \eta(t), \qquad H: \R^{s \times d}, \quad t\ge0,
\end{equation}
where $\eta(t)$ is a squared-integrable function modelling the noise. Consider a filter, that is, the accompanying process described by the following ODE:
\begin{equation} \label{filter}
	\dot{\est} = f(t,\est) + L(t,\est) (y(t) - H(t) \est), \qquad \est(0) = \est_0\,.
\end{equation}
Given $f$ (perfect model), and possibly incomplete observations $y$ ($s<d$), \textbf{the problem is:}
\emph{to find conditions on $H$ which guarantee that there exists a gain $L(t,\est) \in\R^{d \times s}$ such that:}
\begin{equation}
  \label{eq:err_exp_decay}
\|\err(t)\| \le C e^{-a t} \|\err(0)\|\,,\qquad \err\triangleq\truth - \est\,, \qquad t>0\,, C>0\,, a>0\,,
\end{equation}
and \emph{to construct the gain $L(t,\est)$}.

In this work we solve the above problem for generic $f$ by combining ideas from modern Lyapunov stability theory~\cite{BaPe02}, optimal control and numerical analysis. Namely, assuming that the Lyapunov exponents (LEs) associated with the equation $\dot X = \D{f}(t,\est(t))X$ are \emph{forward regular} (see~\cite{BaPe02}), we reformulate the classical notion from control theory, \emph{detectability} for the pair of matrix-valued functions $\matpair{\D{f}(t,\est(t))}{H(t)}$ in terms of LEs, and prove that~\cref{eq:err_exp_decay} holds for every $\err:\|\err(0)\|\le \varepsilon$ if and only if the pair of matrices $\matpair{\D{f}(t,\est(t))}{H(t)}$ is detectable, and $\eta$ disappears when $t\to\infty$. This criterion represents our key contribution from the theoretical standpoint. We then apply this rather theoretical result to construct a gain $L$, and design a numerical algorithm for computing the filter~\cref{filter}. The computation of the gain $L$ relies upon numerically stable procedures for computing LEs of the fundamental matrix differential equation $\dot X = \D{f}(t,\est(t))X$ (see~\cite{DiVV05,DEVV10,DEVV11,DiJoVV11}). Specifically, using a lemma of Perron (cf.~\cite{BaPe02}) and QR-decomposition we design the gain $L$ to guarantee negativity of the LEs associated with the equation governing the estimation error $\err$ in the absence of noise:
\begin{equation}\label{tandyn}
	\dot{\err} = (\D{f}(t,\est(t)) - L(t,\est(t))H(t)) \err + N(t,\err(t)), \quad \err(0) = x_0-\truth_0,\quad\eta\equiv0.
\end{equation}
Negativity of all LEs guarantees \cref{eq:err_exp_decay} under a mild condition on $N$ and for $\|\err(0)\|<\varepsilon$. In other words, the trivial solution of~\cref{tandyn} attracts an $\varepsilon$-neighborhood of $0$.

Instead of solving the ill-conditioned fundamental matrix differential equation directly we compute an orthonormal basis for the non-stable tangent space by solving a matrix differential equation on a manifold of orthogonal matrices~\cite{DiJoVV11}. Since the number of the nonnegative Lyapunov exponents is typically much smaller then the dimension of $\est$, one can reduce the associated computational cost significantly. We stress that this reduction preserves the exponential decay~\cref{eq:err_exp_decay}, and hence the estimation quality is not compromised. The resulting numerical algorithm for computing $L$ and $\est$ is our main contribution from the computational standpoint. Note that our analysis applies to the case where the observations are noise-free ($\eta\equiv 0$). Nevertheless, the proposed filter proves to be efficient in the presence of the observational noise as suggested by our numerical experiments with chaotic non-linear systems.

\paragraph{Motivation and related work}

Problems like \cref{model}--\cref{filter} are fundamental in diverse fields including synchronization in complex networks, data assimilation and control engineering. Syncronizing agents interacting in complex network topologies is of key interest in physical, biological, chemical and social systems to name just a few. In the literature on synchronization of chaotic systems  \cite{PeCa90,PeCa91,PeCaJoMaHe97,BrRu97a,BoKuOsVaZh02}, the coupled system \cref{model}--\cref{filter} is referred to as a driver-receiver process, and the gain $L$ is to be selected so that the receiver $\est$ is synchronized with the driver $\truth$. In their early papers on synchronization, Pecora \& Carroll \cite{PeCa90,PeCa91} note that a necessary condition for synchronization (independent of $x_0$) is that the conditional Lyapunov exponents of \cref{filter} be negative. Similarly to our detectability conditions, conditions for stability of the synchronization manifold, based on a master stability function (MSF) \cite{PeCa98}, also require that all the LEs of an equation similar to~\cref{tandyn} be negative: given an $L$ of the form $L=\nu L_0$ for $\nu$ ranging over the spectrum of the coupling matrix, the MSF assigns a maximal LE of~\cref{tandyn} to each $\nu$. The synchronization manifold is stable provided the MSF maps this spectrum into $(-\infty,0)$. In this work, we use a similar idea to describe a class of observation matrices $H$ for which~\cref{tandyn} has negative LEs: as per~\cref{d:detect-At} below, our detectability condition requires that the non-stable tangent space of $\dot X = \D{f}(t,\est(t))X$ have only trivial intersection with $\ker H^\intercal(t)H(t)$ most of the time. The latter is crucial for the design of the gain $L$, and allows us to achieve~\cref{eq:err_exp_decay}. Note that our algorithm does not require evaluating the MSF for different, possibly complex $\lambda$; instead we need to know the dimension of the unstable tangent subspace. We refer the reader to~\cite{AAGAKJMY_PhRep2008} for a review of recent results on synchronization and the MSF.

In the control literature, the problem of this paper is known as a filtering problem (if $\eta$ and $\truth_0$ are stochastic) or an observer design problem or state estimation problem (if $\eta$ and $\truth_0$ are deterministic). Theoretically, solution of the stochastic filtering problem for Markov diffusions is given by the so-called Kushner-Stratonovich (KS) equation, a stochastic Partial Differential Equation (PDE) which describes evolution of the conditional density of the states of the underlying diffusion process ~\cite{GihmanSkorokhod1997}. For linear systems, the KS equation is equivalent to the Kalman-Bucy Filter~\cite{GhCoTaBuIs81}. In contrast, deterministic state estimators assume that errors have bounded energy and belong to a given bounding set. The state estimate is then defined as a minimax center of the reachability set, a set of all states of the physical model which are reachable from the given set of initial conditions and are compatible with observations. Dynamics of the minimax center is described by a minimax filter. The latter may be constructed by using dynamic programming, i.e., the set $V\le 1$, where $V$ is the so-called value function $V$ solving a Hamilton-Jacobi-Bellman (HJB) equation~\cite{Bardi1997}, coincides with the reachability set~\cite{Baras1995}. Statistically, the uncertainty description in the form of a bounding set represents the case of uniformly distributed bounded errors in contrast to stochastic filtering, where all the errors are usually assumed to be in the form of ``white noise''. However, in many cases $\exp\{-V\}$ coincides with the solution of the KS equation: in fact, for linear dynamics, equations of the minimax filter coincide with those of Kalman-Bucy filter~\cite{KrenerMinimax}.

For generic nonlinear models both minimax and stochastic filters are infinite-dimensional, i.e., to get an optimal estimate one needs to solve either the KS or HJB equation. Hence, if the state space of the model~\cref{model} is of high dimension then both filters become computationally intractable due to the ``curse of dimensionality''. To compute the filter one usually compromises optimality to gain computationally tractable approximations. One such approximation, the Luenberger observer is obtained when the gain $L$ in ~\cref{filter} is chosen so that the estimation error \cref{tandyn} converges to $0$ exponentially, provided $\eta=0$. In fact, there is a deep relationship between observers and (optimal) minimax filters: in the linear case the minimax/Kalman filter uniformly converges to the observer if the observational noise/model error ``disappears'' as $t\to\infty$; see~\cite{BJABJM_SIAMJAM1988}. Motivated by this relationship, we construct an observer $\est$ for \cref{model} as an approximation of the minimax filter. Note that the proposed detectability condition is most related to~\cite{BJABJM_SIAMJAM1988} where a so-called uniform detectability condition\footnote{Uniform detectability: $q^\intercal (\D{f}(t,x) +\Lambda(x) H(t)) q\le -\alpha_0 \|q\|^2$, for some $\alpha_0>0$, $\Lambda$ and all $q,x\in\Rs$.} was used together with uniform controllability to establish~\cref{eq:err_exp_decay}. Similar conditions based on relative degrees are used to design so called high-gain observers~\cite{Deza1992}, minimax sliding mode controllers~\cite{SZAP_TAC16} and minimax filters for differential-algebraic equations~\cite{SZMP_Automatica17,SZMP_TAC16,ZhukSISC13}. A somewhat less restrictive condition of global convergence for skew-symmetric bilinear systems, e.g.~Fourier-Galerkin discretization of the 2D Navier-Stokes equations, was reported in~\cite{SZTTJF_ACC17}. The result of this paper is a lot less restrictive\footnote{$q^\intercal (\D{f}(t,\est(t)) +\Lambda(\est(t)) H(t)) q\ge0$ for all $t$ within any \emph{finite} number of compact intervals of $\R^+$.} than those in the aforementioned papers, as our detectability condition requires that the non-stable tangent space is regularly observed as $t\to\infty$; see~\cref{r:detect}. Note that filters of the form \cref{filter} in the context of Navier-Stokes equations were studied in~\cite{AzOlTi14,GeOlTi16,FoMoTi16} but in the infinite-dimensional setting.

Data Assimilation (DA) improves the accuracy of forecasts provided by physical models and evaluates their reliability by optimally combining \emph{a priori} knowledge encoded in equations of mathematical physics with \emph{a posteriori} information in the form of sensor data. Mathematically, many DA methods rely upon various approximations of stochastic filters. We refer the reader to~\cite{ReichDA2015,StuartDA2015} for further discussions on mathematics behind data assimilation. In what follows we discuss a popular family of algorithms based on Extended Kalman Filter (ExKF) which is the most related to the present work.

ExKF is based on the following idea: given an accurate estimate of the state at time instant $t$, one linearizes the dynamics around that estimate and applies Kalman filtering for the resulting linear system to obtain an estimate for the next time step. This procedure is then repeated. The major drawback of ExKF is that it may diverge for nonlinear equations with positive Lyapunov exponents. A computational bottleneck associated with ExKF is the requirement to recompute the state error covariance matrix. The Ensemble Kalman Filter (EnKF)~\cite{EvensenOceanDyn2003} suggests to overcome this issue by generating an ensemble of trajectories and by computing the ensemble variance to approximate the state error covariance matrix. Recently, this approximation scheme has been combined with ideas from Lyapunov stability theory to construct convergent square-root implementations of the EnKF, a so called EKF-AUS filter~\cite{PaCaTr13,CaGhTrUb08}: the key idea is to sample and propagate the ensemble state error covariance matrix in the unstable tangent subspace only. Apart from solving the ensemble ``deflation problem'', and as noted above, this may provide significant dimension reduction; see~\cite{CaTrDeTaUb08,TrDiTa10}. The importance of Lyapunov exponents and Lyapunov vectors for analysis of DA methods has also been stressed in~\cite{ToHu13}. Very recently, robusteness of the LEs to model errors were studied in~\cite{ColinArxiv} for discrete-time systems.

Finally, we note that computation of LEs is technically challenging~\cite{DiVV05}, as one needs to rely upon some regularity assumptions, e.g.~the aforementioned forward regularity~\cite{BaPe02} or exponential dichotomy~\cite{DEVV11}, to compute LEs for continuos time ODEs. These assumptions are not easy to verify in practice. On the other hand, the regularity for generic discrete time dynamical systems is provided by the so-called Oseledec theorem~\cite{Os68}, hence ``most of the time'' the discrete flow maps resulting from the time/spatial discretizations should be regular. In particular, our experiments with Lorenz '96 (L96) and
Burgers(-Hopf) equations confirmed exponential convergence of the filter in the case of incomplete and noisy observations. We stress that Burgers equation is not dissipative in contrast to the $L96$ system, and it is unclear if this system possesses an invariant ergodic measure making the corresponding discrete flow map subject to the conditions of the Oseledec theorem, so our theory may not apply to this test case. Nevertheless, the proposed filter demonstrates convergence.

\textbf{Notation.} The notation throughout the paper is rather standard. $\|\cdot\|$ denotes the Euclidean norm of $\Rs$, $L^2(0,T)$ is the space of square-integrable functions with norm $\|\cdot\|_{L^2(0,T)}$. $\D{f}$ and $\mathrm{D}^2 f$ denote the Jacobian and Hessian of a smooth map $f:\Rs\to \Rs$.\\
This paper is organised as follows. \Cref{sec:revi-lyap-stab} briefly recalls the notions of Lyapunov exponents, Lyapunov vectors and forward regularity, and reviews how one can compute LEs by means of continuous QR-decomposition. \Cref{sec:detect-line-syst} introduces the notion of detectability for linear non-autonomous systems in terms of LEs,  \Cref{sec:detect-gener-nonl} presents the design of the gain $L$ and the filter $\est$. \Cref{sec:examples} illustrates the application of the filter to Lorenz '96 and Burgers equations. Conclusions are in \Cref{sec:concluding-remarks} and \Cref{sec:proofs} contains the proofs.

\section{Review of Lyapunov stability theory}
\label{sec:revi-lyap-stab}
The stability theory of Lyapunov pertains to linear non\-autonomous differential equations
\begin{equation}\label{linODE}
	\dot{\xi} = A(t) \xi, \qquad \xi(t)\in\Rs,
\end{equation}
where $A(t)\in\R^{d\times d}$ is continuous and bounded.  Following the exposition of Barreira \& Pesin \cite{BaPe02}, we define the function
\begin{equation} \label{LEfun}
\lambda(\xi_0) = \limsup_{t\to\infty} \frac{1}{t} \log \| \xi(t;\xi_0) \|\,,
\end{equation}
which measures the asymptotic rate of exponential growth or decay of the solution of \cref{linODE} with initial condition $\xi(0) = \xi_0$.  Considered over all $\xi_0\in\R^d$, this function assumes $s$ distinct values $\tilde\lambda_1,\dots,\tilde\lambda_s$, $s\le d$. To the function $\lambda$ we assign a \emph{filtration} $\mathcal{V}:=\{V_i\}$, where $V_i:=\{v\in\Rs: \lambda(v)\le \tilde\lambda_i\}$. The filtration $\mathcal{V}$ has the following properties: $V_0\subsetneq V_1\subsetneq,\dots,\subsetneq V_s=\Rs$, and $\lambda(v)=\tilde\lambda_i$ provided $v\in V_i\setminus V_{i-1}$. Letting $n_i:=\operatorname{dim}(V_i)$, $i\ge 1$, $n_0=0$, the number $d_i:=n_i-n_{i-1}$, $i\ge1$, is the multiplicity of $\tilde\lambda_i$ . Hence, the function $\lambda$ assumes $d=d_1+\dots+d_s$ values up to multiplicity: the Lyapunov exponents $\lambda_1 \ge \lambda_2 \ge \cdots \ge \lambda_d$. A basis $\{v_i\}_{i=1}^d$ of $\Rs$ is called a \emph{normal Lyapunov basis} if any $V_i$ can be represented as a linear span of $\{v_{i_1},\dots,v_{i_{n_i}}\}$, $n_i:=\operatorname{dim}(V_i)$. The normal basis is ordered if $V_i$ is equal to  the span of $\{v_{1},\dots,v_{n_i}\}$. Clearly, $\lambda(v_j)=\lambda_j$ for all $n_{i}-d_i< j\le n_i$, $i=1,\dots,s$.

The stability theory of Lyapunov states (see~\cite[p.~9]{BaPe02}) that the trivial solution $\xi(t)=0$ of~\cref{linODE} is exponentially stable if and only if $\lambda_1<0$, that is:
\begin{equation}
  \label{eq:exp_decay}
\forall \varepsilon>0\quad\exists C_\varepsilon>0:\qquad \|\xi(t;\xi_0)\|\le C_\varepsilon e^{(\lambda_1+\varepsilon)t}\|\xi_0\|\,, \qquad \forall t\ge 0\,.
\end{equation}

\subsection{Computation of Lyapunov exponents}
\label{sec:comp-lyap-expon}
The fundamental matrix equation associated to \cref{linODE} is
\begin{equation}\label{fundmat}
	\dot{X} = A(t) X(t), \quad X(t) \in \R^{d\times d},
\end{equation}
whose solution for the initial condition $X(0) = I$ yields the exact flow of \cref{linODE}. If (instead) the columns of $X(0) = (x_1(0), \dots,x_d(0))$ form an ordered normal Lyapunov basis, one finds that $\lambda_i = \lambda(x_i(0))$, $i=1,\dots,d$. The vector $x_i(t)$ is referred to as the $i$th Lyapunov vector corresponding to the $i$th Lyapunov exponent.

The fundamental matrix equation \cref{fundmat} is numerically ill-conditioned, and stabilized formulations are used to compute Lyapunov exponents \cite{DiRuVV97,DiJoVV11,DiVV15}. The algorithm used in this paper relies upon Perron's lemma\cite[Lemma~1.3.3]{BaPe02} which suggests to transform~\cref{fundmat} to the upper-triangular form $\dot R=BR$ by means of a Lyapunov coordinate transformation: $R=Q^\intercal X$. The Lyapunov exponents of $\dot R=BR$ coincide with those of~\cref{fundmat}. Moreover, under a regularity assumption, they can be computed by ``averaging'' the diagonal elements of $B$. For the convenience of the reader we sketch out this procedure below.

Recall from~\cite[p.246]{golub2012matrix} that for any $Y\in\R^{d\times k}$ there exists a QR-decomposition, i.e.~an orthogonal matrix $\fatQ\in\R^{d\times d}$ and upper-triangular matrix $\fatR\in\R^{d\times k}$ such that $Y=\fatQ\fatR$. If the columns of $Y$ are linearly-independent, then, by using the modified Gram-Schmidt (mGS) algorithm\footnote{If the columns of $Y$ are nearly linearly-dependent, and so the condition number of $Y$ is large, the mGS algorithm may generate $\fatQ$ which is not quite orthogonal. In this case a more computationally expensive Housholder transformation can be used instead of the relatively cheap mGS~\cite[p.255]{golub2012matrix} or one may apply mGS twice \cite{Hoffmann89,FrVu99}.}, one can also construct a so called \emph{thin QR-decomposition}: $Y=QR$ where $Q\in\R^{d\times k}$ and $R\in\R^{k\times k}$. We stress that the thin QR-decomposition is unique if one chooses $R_{ii}>0$. Moreover, the ``full'' QR-decomposition is related to the thin one by
\begin{equation}  \label{eq:fat-skin-Q}
	\fatQ = \begin{bmatrix} Q & Q_\perp \end{bmatrix}, \quad \fatR =  \begin{bmatrix} R \\ 0 \end{bmatrix}\,, \quad Q_\perp^\intercal Q_\perp = I\,, Q^\intercal Q_\perp = 0\,.
      \end{equation}
where the columns of $Q_{\perp}\in\R^{d\times(d-k)}$ span the orthogonal complement of the range of $Q$ in $\R^{d\times d}$. If $k=d$ and the columns of $Y$ are linearly-independent, then $\fatQ=Q$ and $\fatR=R$.\\
If $Y\in\R^{d\times k}$, $k\le d$, and $s=\Rank(Y)<k$, then a simple modification of the Gram-Schmidt procedure\footnote{If the columns of $Y^{(k)}:=[Y_1,\dots,Y_k]$ are linearly-independent, and $Y^{(k)}=Q^{(k)}R^{(k)}$ is the thin QR-decomposition of $Y^{(k)}$, and $Y_{k+1}=Q^{(k)}x$, then we set $Q^{(k+1)}:=Q^{(k)}$ and $R^{(k+1)}:=[R^{(k)} x]$. If $Y_{k+2}$ is again in the range of $Q^{(k+1)}$ we repeat the above steps, otherwise we add a row of zeros to the bottom of $R^{(k+1)}$ and append a column of size $k+1$ on the right, and $Q^{(k+2)}$ is computed from $[Y^{(k)} Y_{k+2}]$ by using the mGS procedure.} generates an orthogonal matrix $Q\in\R^{d\times s}$ and an upper-triangular matrix $R\in\R^{s\times s}$ such that $Y=QR$. In this case the diagonal of $R$ may contain zeros.

Now, let $X\in\R^{d\times d}$ solve~\cref{fundmat}, and consider its QR-decomposition $X(t) = \fatQ(t) \fatR(t)$.
Since $\Rank(X(t))\equiv d$, it follows that this QR-decomposition is unique, and $\fatR_{ii}>0$, $i=1,\dots,d$. Moreover, the matrices $\fatQ$ and $\fatR$, constructed by the mGS algorithm, are continuously differentiable. It is easy to find, by differentiating the equality $X = \fatQ\fatR$, that $\fatQ$ and $\fatR$ solve\footnote{Note that, by construction, $\fatR$ is upper-triangular, and so is $\dot \fatR$. This implies that $B$ is upper-triangular too, hence the expression for $S_{ij}$ in~\cref{eq:QR_Q}. } the following system of differential equations:
\begin{align}
    \dot R(t) &= B(t) R(t)\,,\quad R(0) = \fatR_0\,, \quad B=Q^\intercal A Q - S\,,\label{eq:QR_R}\\
    \dot Q(t)& = Q(t) S(t)\,, \quad Q(0) = \fatQ_0\,, \quad S = -S^\intercal\,, S_{ij}= Q_i^\intercal A Q_{j}\,, i>j\,,\label{eq:QR_Q}
\end{align}
where $X(0)=\fatQ_0\fatR_0$. On the other hand, it is not hard to prove that~\cref{eq:QR_R,eq:QR_Q} has a unique solution\footnote{Indeed, \cref{eq:QR_R} is linear in $R$, and $B$ is a continuous function of $t$, hence \cref{eq:QR_R} has a unique solution. The solutions of~\cref{eq:QR_Q} are sought on $O(d)$---the compact manifold of orthogonal square matrices, and since the r.h.s.~of \cref{eq:QR_Q} is Lipschitz in $Q$ w.r.t.~the Frobenius norm on $O(d)$, it follows by the Gr\"{o}nwall-Bellman inequality that there is a unique orthogonal square matrix satisfying~\cref{eq:QR_Q}. Therefore, the unique solution of~\cref{eq:QR_R,eq:QR_Q} coincides with the matrices $\fatQ$ and $\fatR$.} which coincides with the matrices $\fatQ$ and $\fatR$. Hence, the following statements are equivalent:
\begin{enumerate}[label= QR\arabic*, align=left, leftmargin=*]
\item \label{QR1} $X(t)=\fatQ(t)\fatR(t)$ is the unique QR-decomposition of $X$, $\fatQ(t)\fatQ^\intercal(t)=I$, and $\fatR$ is upper-triangular such that $\fatR_{ii}>0$, provided $X$ solves~\cref{fundmat}, and $X(0)=\fatQ_0\fatR_0$.
\item \label{QR2} $\fatQ$ and $\fatR$ solve \cref{eq:QR_R,eq:QR_Q}.
\end{enumerate}
Slightly reformulating Lemma 1.3.5 from~\cite[p.~21]{BaPe02} we introduce the following definition.
\begin{definition}[forward regularity]
Let $\fatR=[R_1,\dots,R_d]$ denote the unique solution of~\cref{eq:QR_R}. The function $r\mapsto\kappa(r_0):=\limsup_{t\to\infty} \frac1t\log(\|R_j(t;r_0)\|)$, where $R_j(t;r_0)$ is the $j$th column of $\fatR(t)$ for $R_j(0)=r_0$, is called \emph{forward regular} provided
\begin{equation}\label{eq:forward_reg}
\limsup_{t\to\infty} \frac1t\int_0^t B_{ii}(s) ds = \liminf_{t\to\infty} \frac1t\int_0^t B_{ii}(s) ds\,,
\end{equation}
\end{definition}
We stress that for bounded $A(t)$, i.e.~$\sup_{t\ge0}\|A(t)\|<+\infty$, Perron's lemma guarantees that $\fatR(t)=\fatQ^\intercal X(t)$ is a Lyapunov transformation which preserves regularity and Lyapunov exponents~\cite[p. 49, Thm. 3.3.1]{Adrianova}, that is:
\emph{$\lambda(\xi_0)$, defined by~\cref{LEfun}, is forward regular, and its range coincides with the range of $\kappa(r_0)$, provided~\cref{eq:forward_reg} holds true.} Moreover, if $\lambda(\xi_0)$ is forward regular then~\cite[p. 24, Thm. 1.3.1]{BaPe02}, for any ordered normal Lyapunov basis $\{v_i\}$ of $\Rs$ it holds that:
\begin{equation}
  \label{eq:LE_QAQ}
\lambda_i = \lambda(v_i) = \lim_{t\to+\infty} \frac 1t \int_0^t B_{ii}(s) ds = \lim_{t\to+\infty} \frac 1t \int_0^t Q_i^\intercal(s) A(s) Q_i(s) ds\,, \;i=1,\dots,d\,.
\end{equation}
Finally, we stress that the orthogonalization process ensures that the Lyapunov exponents are ordered:  $\lambda_1 \ge \cdots \ge \lambda_d$. This important fact allows us to directly compute a basis for the non-stable tangent space, i.e.~the columns $Q_1,\dots,Q_k$ of $\fatQ$ corresponding to $\lambda_1\ge\dots\lambda_k\ge0$, if we know its dimension $k\le d$ a priori, rather than solving \cref{eq:QR_R,eq:QR_Q} for $\fatQ$ and $\fatR$ as suggested above. Indeed, let $X(t) \in \R^{d\times k}$, where $k\le d$, and assume that its thin QR-decomposition is given by:
\[
	X(t)  = Q(t)R(t), \qquad Q(t) \in\R^{d\times k}, \quad R(t)\in\R^{k\times k}.
\]
By differentiating the above equality we find that the equations
\begin{align}
	\dot{Q} &= (I-QQ^\intercal)AQ + Q\tilde S, \label{qr1}\\
	\dot{R} &= \tilde B R, \label{qr2}
\end{align}
can be used to compute $Q\in\R^{d\times k}$ and $R\in\R^{k\times k}$, where as before the skew-symmetric $\tilde{S}^\intercal=-\tilde{S}\in\R^{k\times k}$ is chosen to ensure that $\tilde{B} = Q^\intercal AQ - \tilde{S} \in\R^{k\times k}$ is upper triangular.
If $k=d$, the first term on the right in \cref{qr1} vanishes, and \cref{qr1,qr2} reduce to \cref{eq:QR_R,eq:QR_Q}. The key computational advantage here is that a substantial dimension reduction is achieved when solving~\cref{qr1,qr2} instead of~\cref{eq:QR_R}, provided $k\ll d$.

When $X(0) = Q(0)R(0)$ is of dimension $d\times k$, and assuming that $X(0)$ has nontrivial projection onto the first $k$ elements of some ordered normal Lyapunov basis, the span of the columns of $Q(t)$, found from~\cref{qr1}, will coincide with the span of the Lyapunov vectors corresponding to the leading $k$ Lyapunov exponents. The latter are gleaned from the diagonal of $\tilde B(t)$ via the relation~\cref{eq:LE_QAQ}.

\section{Detectability for linear systems}
\label{sec:detect-line-syst}

In this section we recall the notions of \emph{detectability} and \emph{linear observers} for linear autonomous systems, and we reveal the fundamental relationship between detectability, Lyapunov exponents and existence of linear observers (\cref{l:luenberger}). This relationship is a well-known fact from the classical control theory, and we reformulate it here in terms of the QR-decomposition discussed in~\Cref{sec:comp-lyap-expon} to illustrate using a simple example how to design a linear observer by using Lyapunov vectors and Lyapunov exponents. In particular, we prove (see~\cref{l:luenberger,p:gain}) that the linear observer exists if and only if the system is detectable, and that for non-detectable systems it is impossible, in principle, to design a linear observer as there exists an initial condition such that the corresponding estimation error does not converge to $0$. This simple yet powerful observer design tool is then used in~\Cref{sec:detect-gener-nonl} to construct linear observers for generic non-linear dynamical systems.

\subsection{Autonomous linear systems}
\label{sec:auton-line-syst}

Let us recall the definition of a linear observer, a so-called Luenberger observer. Assume that we observe a function $y(t)\in\Ro$ which represents a possibly incomplete output of a linear system:
\begin{equation}
\label{eq:LTI}
\begin{split}
&\dot \truth(t)  = A\truth(t)\,, \quad \truth(0) = v \in\Rs\,,\\
&y(t) = H\truth(t)\,.
\end{split}
\end{equation}
\begin{definition}\label{d:luenberger}
  The following non-homogeneous linear system
\begin{equation}
  \label{eq:Luenberger}
\dot \est = A\est + L(y - H\est)\,,\quad \est(0) = \est_0\,.
\end{equation}
is called a \emph{linear observer} with gain $L$, the Luenberger gain, provided $L$ is chosen so that the estimation error $\err:=\truth-\est$ converges to $0$ exponentially for any $\est_0$.
\end{definition}
\begin{remark}\label{r:LQreconstruction}
Suppose that $y(t) = H\truth(t)+\eta(t)$ for some $\eta\in L^2(0,T)$, and define $
J(v):=\int_0^T\|y(t) - He^{At}v\|^2\,dt$. Clearly, the best least-squares estimate of the trajectory $\truth$ given $y(t)$ can be computed by minimizing $J$ w.r.t.~$v\in\Rs$ and choosing $z(t) = \exp(At) v$ for this minimizer.
Define $W(0,T)\triangleq\int_0^T e^{A^\intercal t} H^\intercal H e^{At} dt$ and set $b\triangleq\int_0^T e^{A^\intercal t} H^\intercal y(t) dt$. We find that $J(v) = \|y\|^2_{L^2(0,T)} - 2 v^\intercal b + v^\intercal W(0,T)v$, and thus the set of minimizers of $J$ is given by $\{v:W(0,T)v = b\}$. For instance, the minimizer with smallest 2-norm is given by: $v^\dag = W^\dag(0,T)b$. It is not hard to see that $J(v^\dag)= \|\eta\|^2_{L^2(0,T)}+ \widetilde{\eta}^\intercal W^\dag(0,T) \widetilde{\eta} + 2 \widetilde{\eta}^\intercal (I-W^\dag(0,T)W(0,T)) \truth(0)$, where $\widetilde{\eta}\triangleq \int_0^T e^{A^\intercal t} H^\intercal \eta(t) dt$. If $\eta=0$, one can reconstruct $\truth$ exactly if and only if $W(0,T)$ is of full rank: indeed, $J(v^\dag)=0$ in this case. If $\eta=0$ but $W(0,T)$ is not of full rank, any solution $v$ of $W(0,T) v = b$ can be represented as $v=v^\dag\oplus v_0$ where $W(0,T)v_0=0$. We stress that $v_0$ cannot be determined from $W(0,T) v = b$. In fact, $H e^{At}v_0=0$, and so the observed data $y$ does not contain any information about the null-space of $W(0,T)$. In the control literature, this sub-space is sometimes referred to as an unobservable subspace. Finally, if $\eta\ne0$ one cannot reconstruct $\truth$ exactly, independent of the rank of $W$: the best estimate of $\truth$ will have the mean-squared error $J(v^\dag)$.
\end{remark}
The unobservable subspace can be computed efficiently for the case of linear autonomous systems. Define \[O^s:=\left[
  \begin{smallmatrix}
    H\\ HA\\\vdots\\HA^{s}
  \end{smallmatrix}
\right]\,, \quad \text{$s$ the smallest integer such that $\operatorname{rank}(O^{s}) = \operatorname{rank}(O^{s+p})$, $\forall p\ge 1$\,.}
\]
In turns out (see proof of~\cref{l:luenberger}) that $\ker W(0,T)=\ker O^s$. Note that the state space of~\cref{eq:LTI} splits into two invariant sub-spaces: $\Rs=\ker O^s\oplus \ker^\perp O^s$, and that the part of any trajectory~\cref{eq:LTI} in $\ker^\perp O^s$, i.e.~the projection of the state vector $\truth(t)$ onto $\ker^\perp O^s$, can be recovered from the data $y$. The ``invisible'' part, i.e.~the projection of $\truth(t)$ onto $\ker O^s$, cannot be recovered unless the matrices $A$ and $H$ have a specific structure:
\begin{definition}\label{d:detect_LTI}
$\matpair{A}{H}$ is \emph{detectable} if $\lim_{t\to\infty}\|e^{At}v\|=0$ for any $v\in\ker(O^s)$.
\end{definition}
Detectability implies that the projection of $\truth(t)$ onto $\ker^\perp O^s$, the ``invisible'' part, decays to zero exponentially fast. The following lemma establishes the connection between detectability, existence of linear observers and Lyapunov exponents.
\begin{lemma}\label{l:luenberger}
$\matpair{A}{H}$ is detectable if and only if all the Lyapunov exponents corresponding to $\ker(O^s)$ are negative, i.e.~$\lambda(v)<0\,, \forall v\in\ker(O^s)$. If $\matpair{A}{H}$ is detectable then there exists a linear observer~\cref{eq:Luenberger} in the sense of~\cref{d:luenberger}. If, on the contrary, there exist a vector $v^\dag\in \ker(O^s)$ such that the corresponding Lyapunov exponent is non-negative, $\lambda(v^\dag)\ge 0$, then for any gain matrix $L$, the estimation error $\err=\truth-\est$ either grows unbounded, i.e.~$\lim_{t\to+\infty}\|e(t)\|=+\infty$, or stays bounded but can be made arbitrarily large. 
\end{lemma}

\subsection{Non-autonomous linear systems}
\label{sec:non-auton-line-syst}
In this subsection we generalize~\cref{l:luenberger} to the case of time-dependent matrices $A$ and $H$. Specifically, a time-variant version of~\cref{eq:LTI} takes the following form:
\begin{equation}
\label{eq:LTV}
\begin{split}
  \dot \truth(t) &= A(t)\truth(t)\,, \quad \truth(0) = v\,,\\
  y(t)&=H(t)\truth(t).
\end{split}
\end{equation}
In what follows we generalize~\cref{d:detect_LTI} to the case of time-dependent $A$ and $H$, and then this definition is applied to design the linear observer for~\cref{eq:LTV}. Assume that $t \mapsto H(t)\in\R^{s\times d}$, $s\le d$ and $t \mapsto A(t)\in\R^{d \times d}$ are such that
\begin{equation}
  \label{eq:AH}
  \sup_{t\ge 0} \sup_{\|x\|=1}\| A(t)x \| < +\infty\,,\quad \sup_{t\ge 0} \sup_{\|x\|=1}\| H(t)x \| < +\infty\,.
\end{equation}
\begin{definition}\label{d:detect-At}
Assume that~\cref{eq:forward_reg} holds true, and let $Q\in\R^{d\times k}$, $k\le d$ solve~\cref{qr1}. Let $\widetilde{Q}(t)\in\R^{d\times k}$ and $\widetilde{R}(t)\in\R^{k \times k}$ be the thin QR-decomposition of $H^{\intercal}HQ$, i.e.
\[
	\widetilde{Q}(t) \widetilde{R}(t) = H^{\intercal}(t) H(t) Q(t).
\]
We say that $\matpair{A}{H}$ is {\bf detectable in the direction $Q_j$} if
\[
	\limsup_{t\to+\infty} \frac{1}{t} \int_0^t \widetilde{R}_{jj}(s)\,ds > 0.
\]
Furthermore we say that $\matpair{A}{H}$ is {\bf detectable} if it is detectable in any direction $Q_j$ from the non-stable tangent space:
\begin{equation}
  \label{eq:LEQA}
	\lambda_j = \lim_{t\to+\infty} \frac{1}{t} \int_0^t B_{jj}(s)\,ds = \lim_{t\to+\infty} \frac{1}{t} \int_0^t Q_j(s)^{\intercal} A(s) Q_j(s)\,ds \ge 0\,.
      \end{equation}
\end{definition}
\begin{remark}\label{r:detect}
Note that $\widetilde{R}_{jj}\ge 0$ by construction, with equality $\widetilde{R}_{jj}(t)=0$ holding only if $H^\intercal(t)H(t)Q(t)$ is rank deficient\footnote{In fact, $\tilde R_{ii}=0$ if and only if the $i$th column of $H^\intercal(t)H(t)Q(t)$ can be represented as a linear combination of the first $i-1$ columns. In this case, the $i$th column of $\tilde Q$ is set to $0$ by the mGS process.} The latter is the case if and only if the linear sub-space generated by the columns of $Q$, $\{Q(t)x,x\in\R^k\}$ has a non-trivial intersection with $\ker H^\intercal(t)H(t)$: $H^\intercal(t)H(t)Q(t)x = 0$ for some $x\in\R^k$. Recall from~\cref{sec:comp-lyap-expon} that the way we compute $Q$ ensures that the Lyapunov exponents are ordered:  $\lambda_1 \ge \cdots \ge \lambda_{k^\star}$, and the $k^\star$ leading columns of $Q(t)$, $Q_1,\dots,Q_{k^\star}$ correspond to the $k^\star$ leading Lyapunov exponents as in~\cref{eq:LEQA}. Now, if $k^\star$ is the number of non-negative Lyapunov exponents of~\cref{eq:LTV}, then, according to~\cref{d:detect-At}, $\matpair{A(t)}{H(t)}$ is detectable iff the sub-space generated by the first $k^\star$ columns of $Q$, i.e.~the non-stable tangent space, has only trivial intersection with $\ker H^\intercal(t)H(t)$ most of the time, i.e.~the measure of the set $\{0\le s\le t:\widetilde{R}_{jj}(s)>0\}$ grows at least linearly as $t\to+\infty$, outside perhaps a \emph{finite} number of compact intervals of $\R^+$, for $j=1,\dots,k^\star$.

Intuitively, detectability requires that the non-stable tangent space, spanned by $Q_1,\dots,Q_{k^\star}$, is regularly observed as $t\to\infty$. If, however, there is a vector $Q_\ell$ in the non-stable tangent space that is `unseen' by the observation operator $H$ most of the time, i.e.~the measure of the set $\{0\le \tau\le t: R_{\ell\ell}(\tau)>0\}$ is finite or grows at a sub-linear rate as $t\to+\infty$, then $\matpair{A(t)}{H(t)}$ is not detectable.
\end{remark}

\cref{r:detect} suggests that \emph{a necessary condition for detectability of $\matpair{A(t)}{H(t)}$} is
\begin{equation}
  \label{eq:Detect_nesess}
  \min\{\Rank(H(t)),k\} \ge k^\star\,, \text{ for almost all } t\in (\R^+\setminus K)\,,
\end{equation}
provided $K$ is a union of a finite number of compact intervals of $\R^+$, and $k^\star$ is the number of nonnegative Lyapunov exponents of \cref{eq:LTV}. Indeed, $\ker H^\intercal(t)H(t) = \ker H(t)$ as the range of $H$ is orthogonal to $\ker H^\intercal$. Hence, $\Rank H^\intercal HQ = \Rank HQ$. By construction, $Q(t)$ has full column rank $k$. Hence, $\Rank H^\intercal(t)H(t)Q(t) = \Rank H(t)Q(t) = \min\{\Rank(H(t)),k\}$ at any time $t$ such that $\ker(H(t)Q(t))=\{0\}$. In the latter case, $\Rank\widetilde{R}(t) = \min\{\Rank(H(t)),k\}$. Thus, the diagonal of $\widetilde R(t)$ may have at most $\min \{ \Rank(H(t)) , k \}$ positive elements for any $t>0$. As noted above, detectability of $\matpair{A(t)}{H(t)}$ implies that, for any $j=1,\dots,k^\star$, the measure of the set $\{0\le s\le t:\widetilde{R}_{jj}(s)>0\}$ grows at least linearly as $t\to+\infty$, outside perhaps a \emph{finite} number of compact intervals of $\R^+$, $K$. Hence, there exist $T>0$ such that $\widetilde{R}_{jj}(s)>0 $ for all $j=1,\dots, k^\star$, and almost all $s>T$. But this is possible only if~\cref{eq:Detect_nesess} holds true.  In particular, \cref{eq:Detect_nesess} implies that \emph{the rank of the observation operator $H(t)$ must equal or exceed the number of nonnegative Lyapunov exponents for detectability to hold.}

The following proposition shows that detectability is necessary and sufficient for the existence of the linear observer.
\begin{proposition}\label{p:gain}
Assume that $\eta=0$, $A$ and $H$ satisfy~\cref{eq:AH}, and \cref{eq:forward_reg} holds true. Let $Q\in\R^{d\times k}$ and $R\in \R^{k\times k}$ solve~\cref{qr1,qr2}, provided $k$ is the number of non-negative Lyapunov exponents of~\cref{fundmat}, and let $\widetilde Q\widetilde R = H^\intercal H Q$ be the thin QR-decomposition of $H^\intercal H Q(t)$. Define the gain
\begin{equation}
  \label{eq:gainL}
  	L := p\, Q\, \widetilde{Q}^{\intercal}  H^\intercal\,, \quad p>0\,,\quad H^\intercal H Q = \widetilde{Q}\widetilde R\,,
\end{equation}
and let $\est$ solve the following system:
\begin{equation}
\label{eq:filter_lin}
  \begin{split}
    	\dot{\est} &= A(t)\est + L(t) (y(t) - H(t) \est(t)), \quad \est(t) = 0\,.
  \end{split}
\end{equation}
Then there exists $p>0$ such that the estimation error $\err=\truth-\est$ decays to $0$ exponentially fast if and only if $\matpair{A}{H}$ is detectable in the sense of~\cref{d:detect-At}.
\end{proposition}
\begin{remark}
  Note that the gain $L$ is designed so that all the non-negative LEs of $A(t)-L(t)H(t)$ are made negative. In other words, the gain $L$ (asymptotically) stabilizes the error equation $\dot e = (A(t) - L(t) H(t)) e$, provided the observations are exact. The latter is possible iff $\matpair{A(t)}{H(t)}$ is detectable as per~\cref{d:detect-At}.
\end{remark}
The proof of \cref{p:gain} is given in the appendix.

\section{Approximation of the minimax filter for nonlinear ODEs}
\label{sec:detect-gener-nonl}


As noted in the introduction, for linear systems the minimax/Kalman filter uniformly (in time) converges to the linear observer, provided the observational noise/model error ``disappears'' as $t\to\infty$; see~\cite{BJABJM_SIAMJAM1988}. In the general case of non-linear dynamics and only bounded observational noise, linear observers represent an approximation of the minimax filter. In this section we design a linear observer $\est$ for a nonlinear ODE
\begin{equation}\label{eq:state}
	\dot{\truth}  = f(t,\truth), \quad \truth(t) \in \mathcal{D}\subset \Rs, \quad \truth(0) = \truth_0\,,
\end{equation}
with compact state space $\mathcal{D}$, given incomplete and noisy\footnote{$t\mapsto \eta(t)\in\Ro$ is a bounded square-integrable function such that $\|\eta(t)\|\le\varepsilon$} observations:
\begin{equation}\label{eq:obs}
	y(t) = H(t) \truth(t) +\eta(t), \quad H: \R^{s \times d}\,, s\le d\,,\quad \eta\in L^2(0,+\infty)\,, \quad t\ge0\,.
\end{equation}
Assuming that~\cref{eq:state} has the unique solution for any $\truth_0\in\mathcal{D}$ we prove that, locally the estimation error $\err=\truth-\est$ decays to $0$ exponentially iff there exists a gain $L(t,\est)$ such that all the LEs of $A-LH$ are negative, provided $A(t):=\D{f}(\est(t))$, and $\eta=0$. If the function $\lambda$ defined by~\eqref{LEfun} is, in addition, forward regular, then the LEs of $A-LH$ are negative iff $A,H$ is detectable, and the gain $L$ can be constructed as per the recipe of~\cref{p:gain}. We study the case of non-trivial noise $\eta$ numerically in the following section.

The following proposition utilizes Lyapunov stability theorem~\cite[p. 29, T.1.4.3]{BaPe02} to prove that under mild regularity conditions a vicinity of $0$ is attracted to $\err\equiv 0$ by the error dynamics. 
Define $N(t,\xi,x):=(N_1(t,\xi,x)\dots N_d(t,\xi,x))^\intercal$ where \[
N_i(t,\xi,x):=\xi^\intercal \left(\int_0^1 \int_0^1 s \mathrm{D}^2 f_i (x + \tau s \xi, t) ds\,d\tau\right)\xi\,.
\]
We require the following assumption on the linearized dynamics of the observer.
\begin{assumption}\label{a:lindyn}
The following conditions hold for the Jacobian $\D f(t,\est(t))$ evaluated along the observer process:
\begin{enumerate}[label= A\arabic*, align=left, leftmargin=*]
\item \label{a1} $A(t):=\D{f}(t, \est(t))$ is bounded: $\sup_{t\ge 0} \sup_{\|x\|=1}\| A(t)x \|< +\infty$.
\item \label{a2} $\max_{x\in K}\|D^2 f_i (x, t)\|\le C_i(K)<+\infty$ for every compact subset $K$ of $\Rs$.
\item \label{a3}the Lyapunov exponent $\lambda(\xi)$ of $\dot X = A(t)X$ is forward regular.
\end{enumerate}
\end{assumption}

\begin{theorem}\label{t:tsf}
Let $\eta\equiv0$ and assume the conditions of \cref{a:lindyn} hold.  Let $\est$ be the unique solution of the following system:
\begin{equation}
\label{eq:filter}
  \begin{split}
    	\dot{\est} &= f(t,\est) + L(t,\est) (y(t) - H(t) \est(t)), \qquad \est(0) = \est_0,
\end{split}
\end{equation}
Then the following statements are equivalent:
\begin{enumerate}[label= S\arabic*, align=left, leftmargin=*]
\item \label{s1} There exist $\varepsilon>0$ and $a>0$ such that $\|\err(t)\|\le C e^{-at}\|\err(0)\|$ for all $t\ge0$ where $\err=\truth-\est$ is the unique solution of the error equation:
  \begin{equation}
    \label{eq:error}
    \dot \err = (A(t) - L(t,\est)H^\intercal(t) H(t))\err + N(t,\err,\est(t))\,, \quad \err(0)=\err_0\,, \|\err_0\|<\varepsilon.
  \end{equation}
\item \label{s2} All Lyapunov exponents of $\dot X = (A-LH)X$ are negative.
\end{enumerate}
\end{theorem}
\begin{remark}
  In fact, \cref{t:tsf} establishes robustness of the LEs w.r.t.~small perturbation. Indeed, one may interpret $\err(t)$ as a perturbation about a particular solution $\est$, i.e.~$\truth(t) = \est + \err$, and the dynamics of this perturbation is defined by the  ``perturbed'' equation~\cref{eq:error}: here $N$ is considered as a perturbation of the linear unperturbed equation $X=(A-LH)X$. By \cref{t:tsf} it follows that the exponential decay of solutions of the ``perturbed'' equation~\cref{eq:error} is equivalent to the exponential decay of the solutions of unperturbed linear equation provided the LEs are forward-regular. Osedelec's theorem \cite{Os68} establishes regularity for a wide class of nonlinear systems possessing an ergodic invariant measure. For a nonlinear system with a global ergodic attractor, the LEs are independent of any particular trajectory $\est$.
\end{remark}
\begin{corollary}\label{c:tsf}
Assume that all the conditions of \cref{t:tsf} hold and let $\est$ solve
\begin{align}
    \dot{\est} &= f(t,\est) + L(t,\est) (y(t) - H(t) \est(t))\,, & L &= p\, Q \widetilde{Q}^{\intercal}  H^\intercal,p>0,\quad \est(0) = \est_0 \label{eq:TSFilter}\\
    \dot{Q} &= (I-QQ^\intercal)\D{f}(t,\est)Q + QS\,, & Q&(0)=Q_0\in\R^{d\times k}\,,\\
 S &= -S^\intercal, S_{ij} = (Q^\intercal \D f(t,x) Q)_{ij},i>j, & H^\intercal& H Q = \widetilde{Q}\widetilde R
\end{align}
for appropriate $k \in 1, \dots, d$. Then the estimation error $\err=\truth-\est$ converges to zero exponentially if and only if $\matpair{A}{H}$ is detectable in the sense of~\cref{d:detect-At}.
\end{corollary}

\section{Numerical experiments with noisy observations}
\label{sec:examples}
In this section we apply~\cref{eq:TSFilter} to the Lorenz '96 (L96) model and Burgers equation, and compare it to the Extended Kalman(-Bucy) Filter (ExKF). We consider both exact and noisy observations. Abusing the standard control terminology we will refer to~\cref{eq:TSFilter} as a filter to stress that in the experiments the observations are allowed to be noisy.

\subsection{L96}
\label{sec:L96}
\begin{figure}[h]
\begin{centering}
\includegraphics[width=0.49\textwidth]{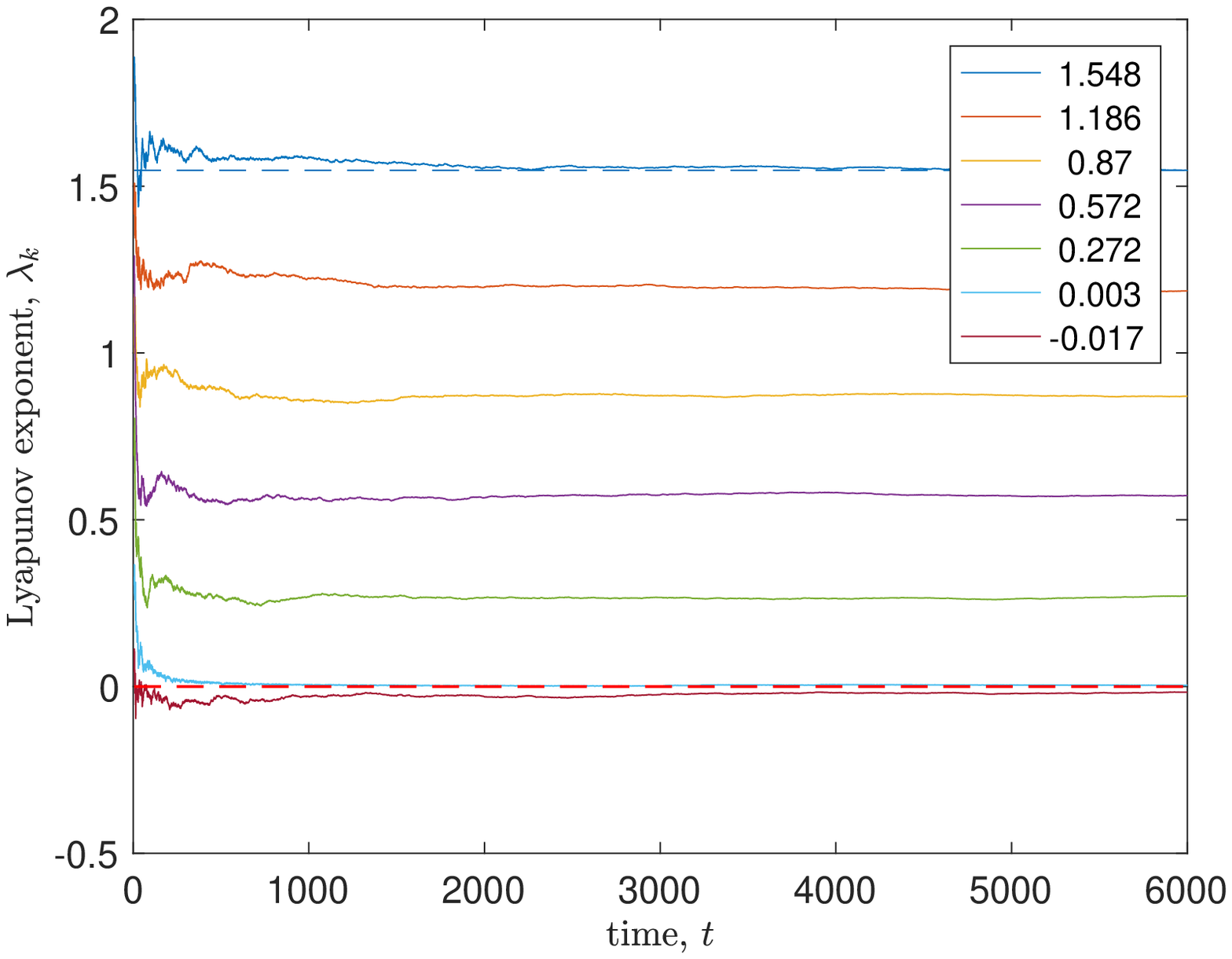}
\includegraphics[width=0.49\textwidth]{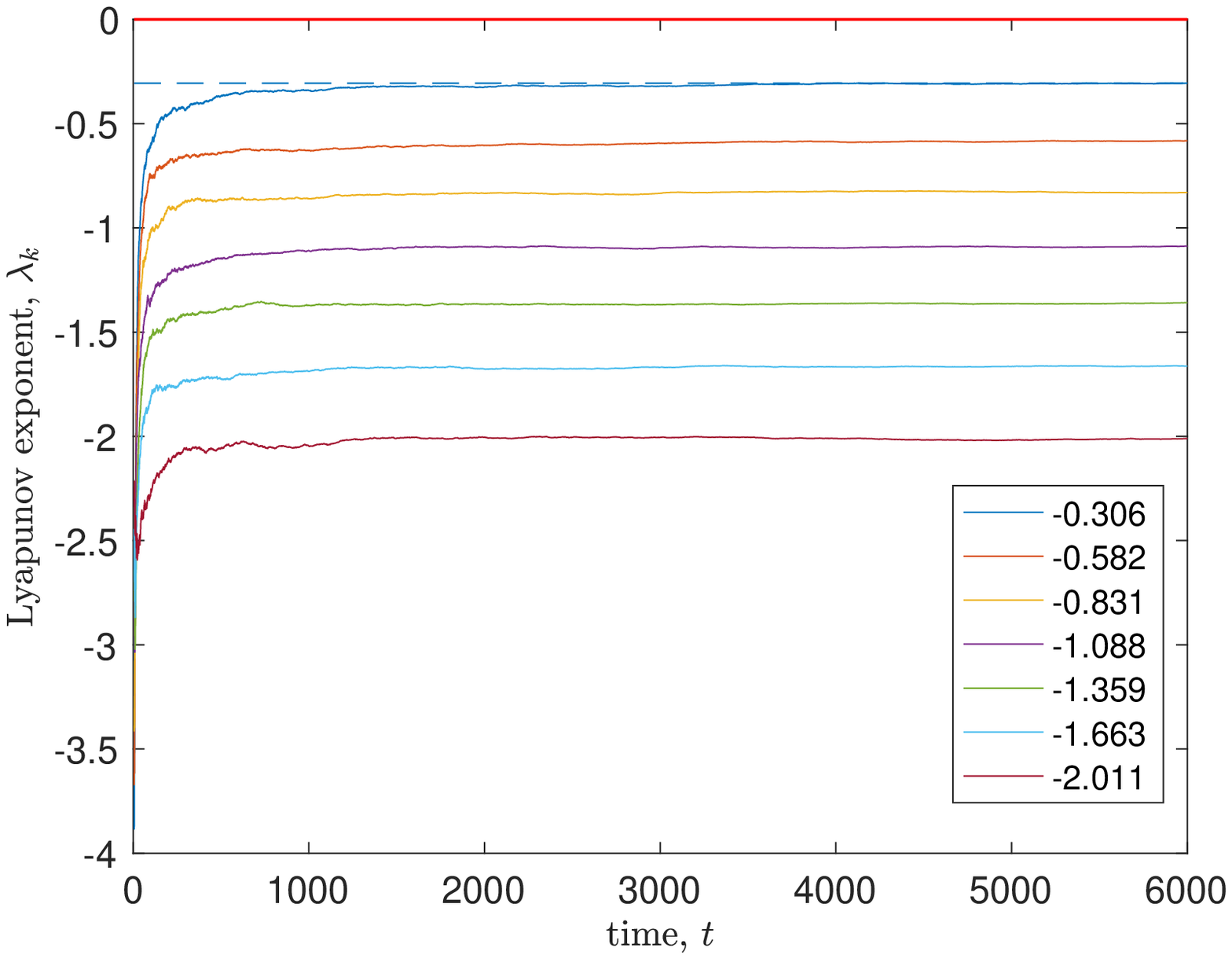}
\caption{Leading Lyapunov exponents of the Lorenz '96 model \cref{eq:L96} (left) and the tangent space dynamics \cref{eq:error} (right).\label{fig3-4}}
\end{centering}
\end{figure}

The Lorenz '96 model is a system of ODEs
\begin{equation}
  \label{eq:L96}
\dot \truth_i = -\truth_{i-2} \truth_{i-1} + \truth_{i-1} \truth_{i+1} -\truth_i +\mathcal{F}\,, \quad \truth_i(0) = \sin\left(2\pi \frac{i-1}{d}\right),~i=1,\dots,d,
\end{equation}
defined on a periodic lattice with state space dimension $d=18$, and constant forcing $ \mathcal{F}\equiv 8$. We solve~\cref{eq:L96} using a fourth order explicit Runge-Kutta method (RK4) with the time step $\Delta t=0.01$ to sample the observations $y$. The rows of the matrix $H$ are taken to be the first $k$ eigen-vectors of the discrete Laplacian (i.e.~$d\times d$ circulant matrix having tridiagonal elements $\left[ \begin{smallmatrix} 1 & -2 & 1 \end{smallmatrix}\right]$). As per estimate~\eqref{eq:Detect_nesess}, the rank of the observation operator $H$, $k$ must be greater or equal to the number of nonnegative LEs of $\dot X = \D{f}(\est(t))X$ along the trajectory $\est$ of~\cref{eq:TSFilter}. We computed RK4 approximations of these LEs by averaging over the interval $t\in[0,6000]$: on the left panel of \cref{fig3-4} one can see that $6$ leading exponents are nonnegative, but the $7$th exponent, $\lambda_7 \approx -0.017$. In addition, for $k=6$ the approximation of the leading exponent of the error dynamics~\eqref{eq:error} equals $-0.051$. The latter two observations suggest that the choice $k=6$ may lead to very slow convergence or convergence only within a very small neighborhood of the truth $x(t)$, and, for the ``boundary'' case $k=6$, the convergence is very sensitive both to the size of the initial perturbation, and to the accuracy with which one is able to approximate the basis $Q(t)$ of the non-stable tangent space. To alleviate this sensitivity, we set the rank of $H$ to $k=7$. The latter ensures that indeed, as per condition \ref{s2} of~\cref{t:tsf}, all the exponents of the error dynamics~\cref{eq:error} are negative, and, more importantly, they are well separated from $0$; see the right panel of \cref{fig3-4}.

\begin{figure}[h]
\begin{centering}
\includegraphics[width=0.52\textwidth]{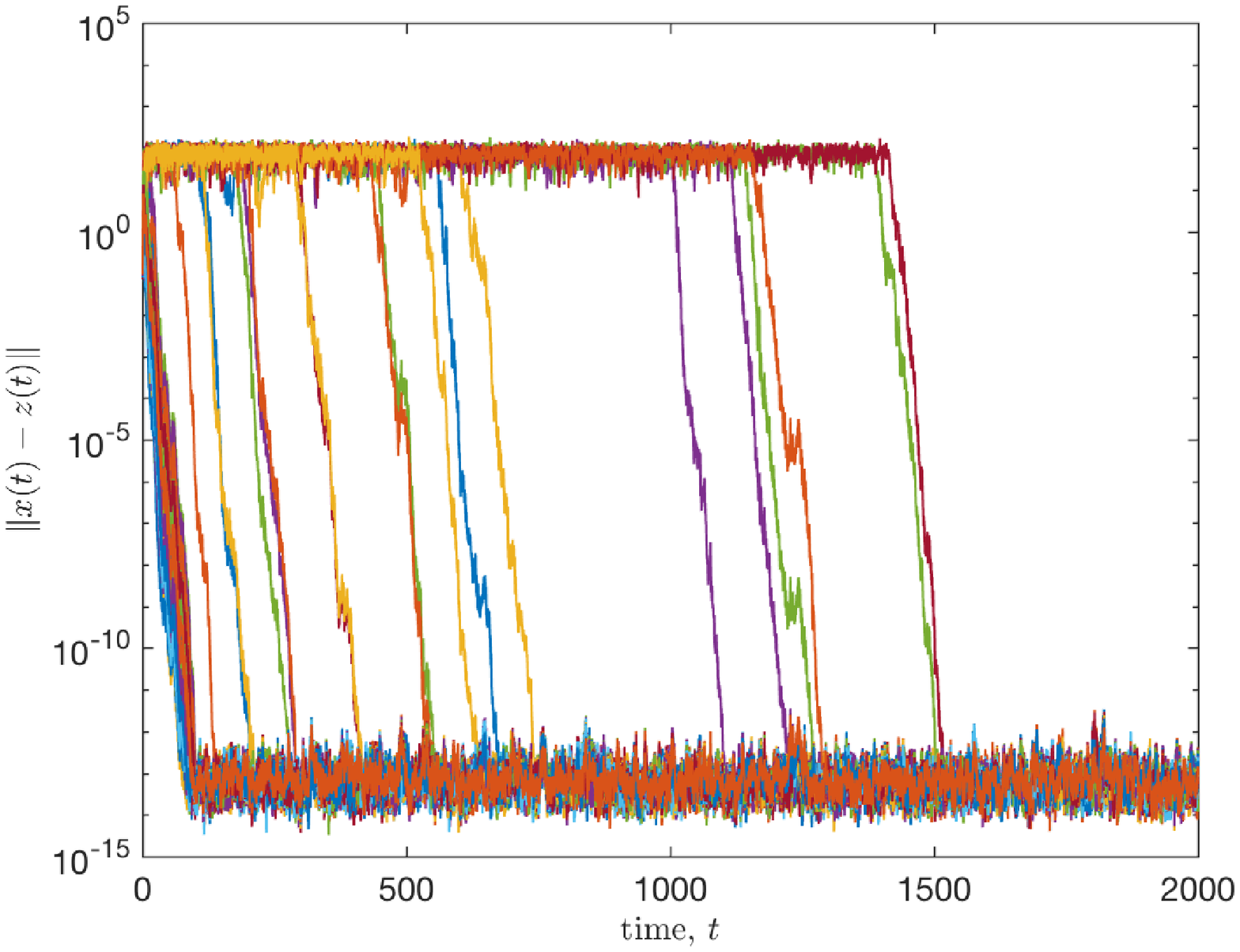}
\includegraphics[width=0.47\textwidth]{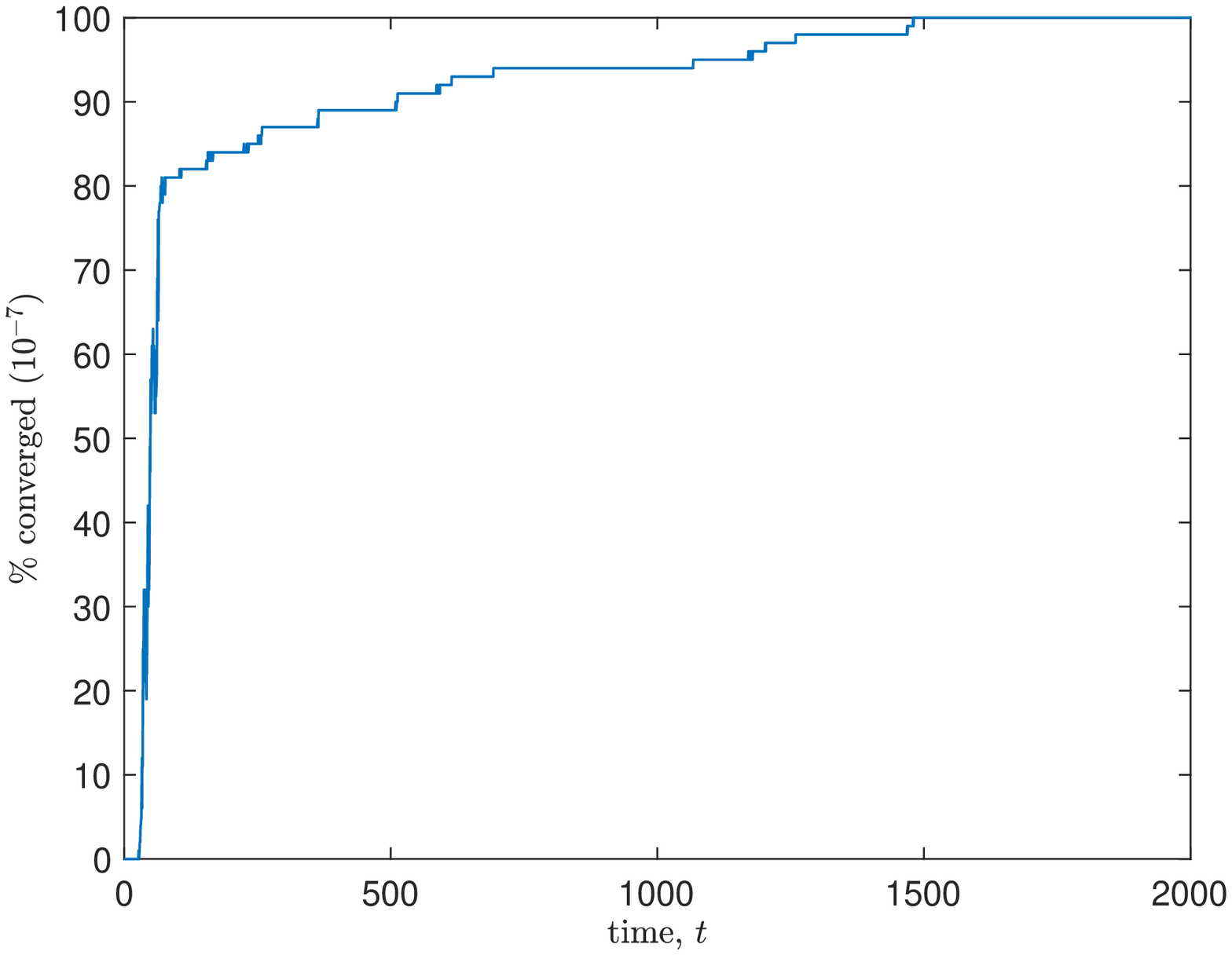}
\caption{Convergence of filter \cref{eq:TSFilter} for the Lorenz '96 model \cref{eq:L96} with $k=7$.  Left, the errors $\|\xi(t)\|$ for a 100-member ensemble of perturbed initial conditions.  Right, the number of ensemble members converged to tolerance $\|\xi(t) \|<10^{-7}$ at time $t$.\label{fig1-2}}
\end{centering}
\end{figure}

For the case $k=7$ we generate an ensemble of 100 initial conditions $\est_0=\truth(0)+0.01\eta$, $\eta_i\sim\mathcal{N}(0,1)$, $i=1,\dots,d$. For each $\est_0$ we set $p=10$, and integrate~\cref{eq:TSFilter} using the same time integrator until $\|\est-\truth\|\le 10^{-14}$. The resulting 2-norm estimation error $\| \xi(t) \| = \| x(t)-z(t) \|$ as a function of time are shown in the left panel of \cref{fig1-2}. We see that all samples ultimately converge at an exponential rate, as predicted by the theory of this paper.  However, in some cases there is a long delay before exponential convergence is observed.  In the right panel of \cref{fig1-2} we plot the number of ensemble members that has converged to within tolerance $\|\xi(t) \| < 10^{-7}$ at time $t$.  We see that 80\% of the ensemble converges by time $t=100$, the remaining 20\% converges more slowly, with the last ensemble member converging only at time $t=1500$.

\begin{figure}[h]
\begin{centering}
\includegraphics[width=0.52\textwidth]{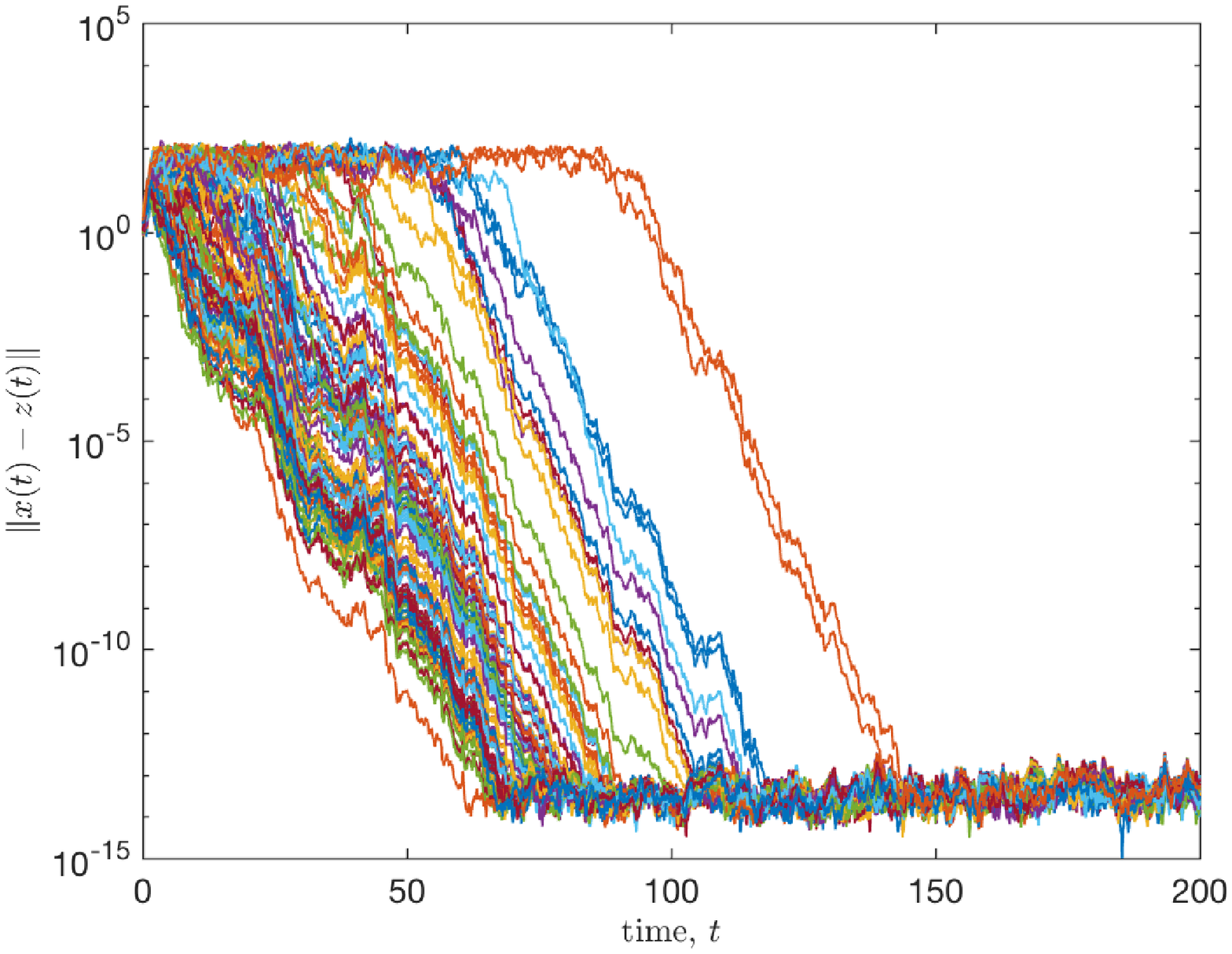}
\includegraphics[width=0.47\textwidth]{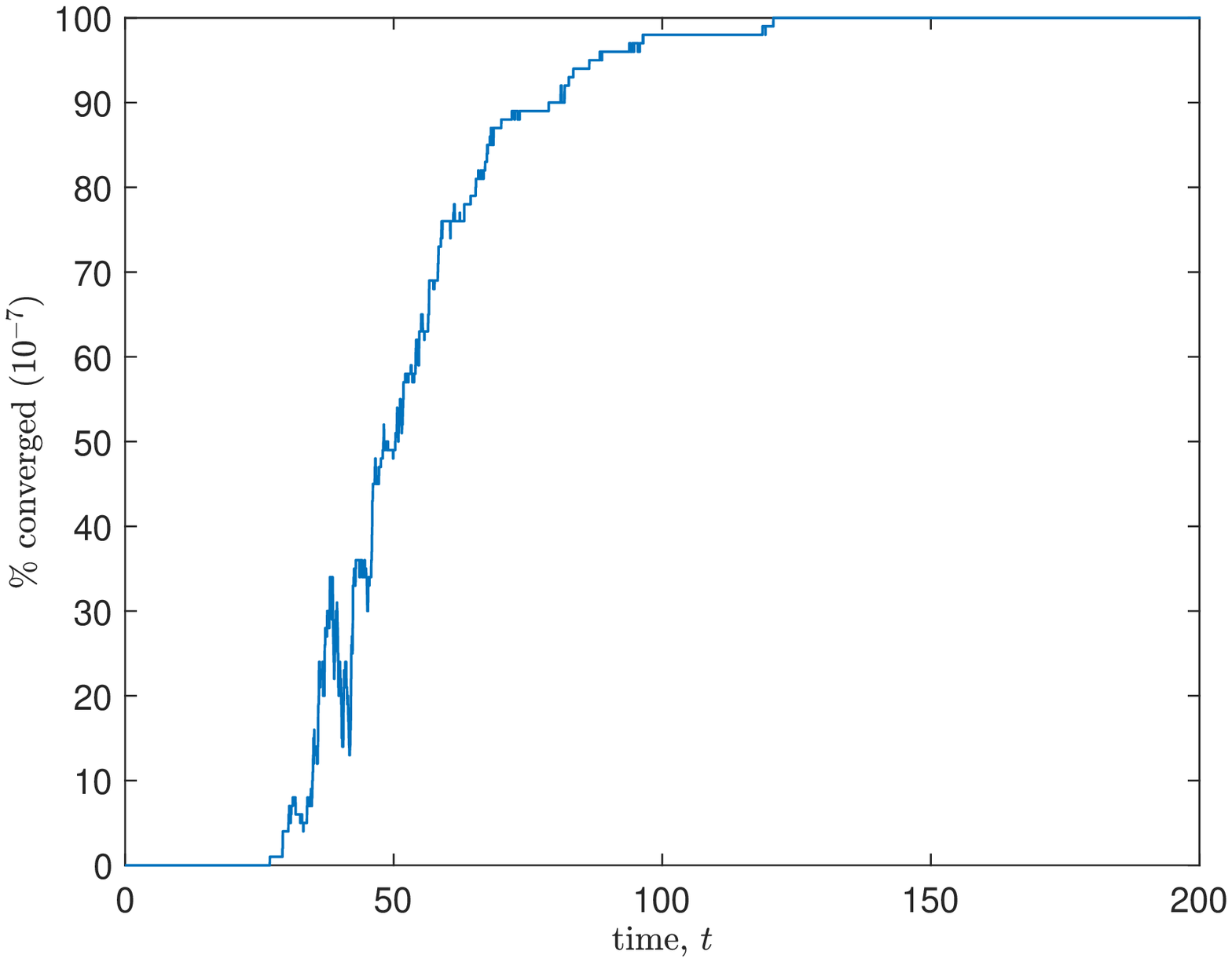}
\caption{Convergence of filter \cref{eq:TSFilter} for the Lorenz '96 model \cref{eq:L96} with $k=8$.  Left, the errors $\|\xi(t)\|$ for a 100-member ensemble of perturbed initial conditions.  Right, the number of ensemble members converged to tolerance $\|\xi(t) \|<10^{-7}$ at time $t$.\label{fig5-6}}
\end{centering}
\end{figure}

To make the convergence of~\cref{eq:TSFilter} faster and even less sensitive to the size of the initial perturbation, and to the errors of RK4 approximation of the basis $Q(t)$, we set $k=8$ and increase the amplitude of the initial perturbation to
$\est_0=\truth(0)+0.1\eta$, $\eta_i\sim\mathcal{N}(0,1)$, $i=1,\dots,d$ so that $\frac{\|\truth_0-\est_0\|}{\truth_0}\approx 0.3$, i.e.~the magnitude of the perturbation is up to $30\%$. Figure \ref{fig5-6} demonstrates the convergence. Even with a much larger initial perturbation, all ensemble members converge to machine precision within time $t=200$, and more than 95\% converge to within $\|\xi(t)\|<10^{-7}$ by time $t=100$.

\begin{figure}[h]
\begin{centering}
\includegraphics[width=0.49\textwidth]{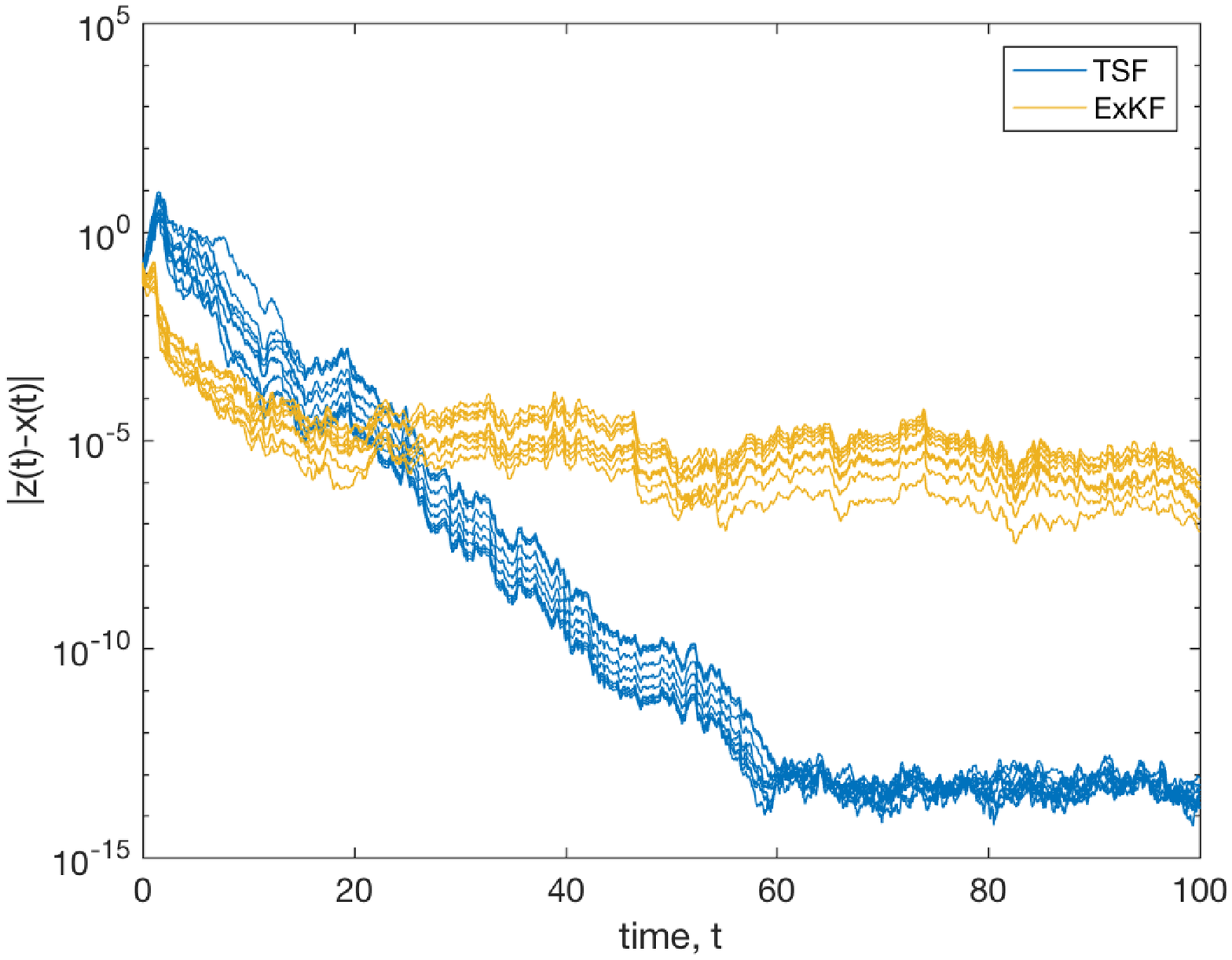}
\includegraphics[width=0.49\textwidth]{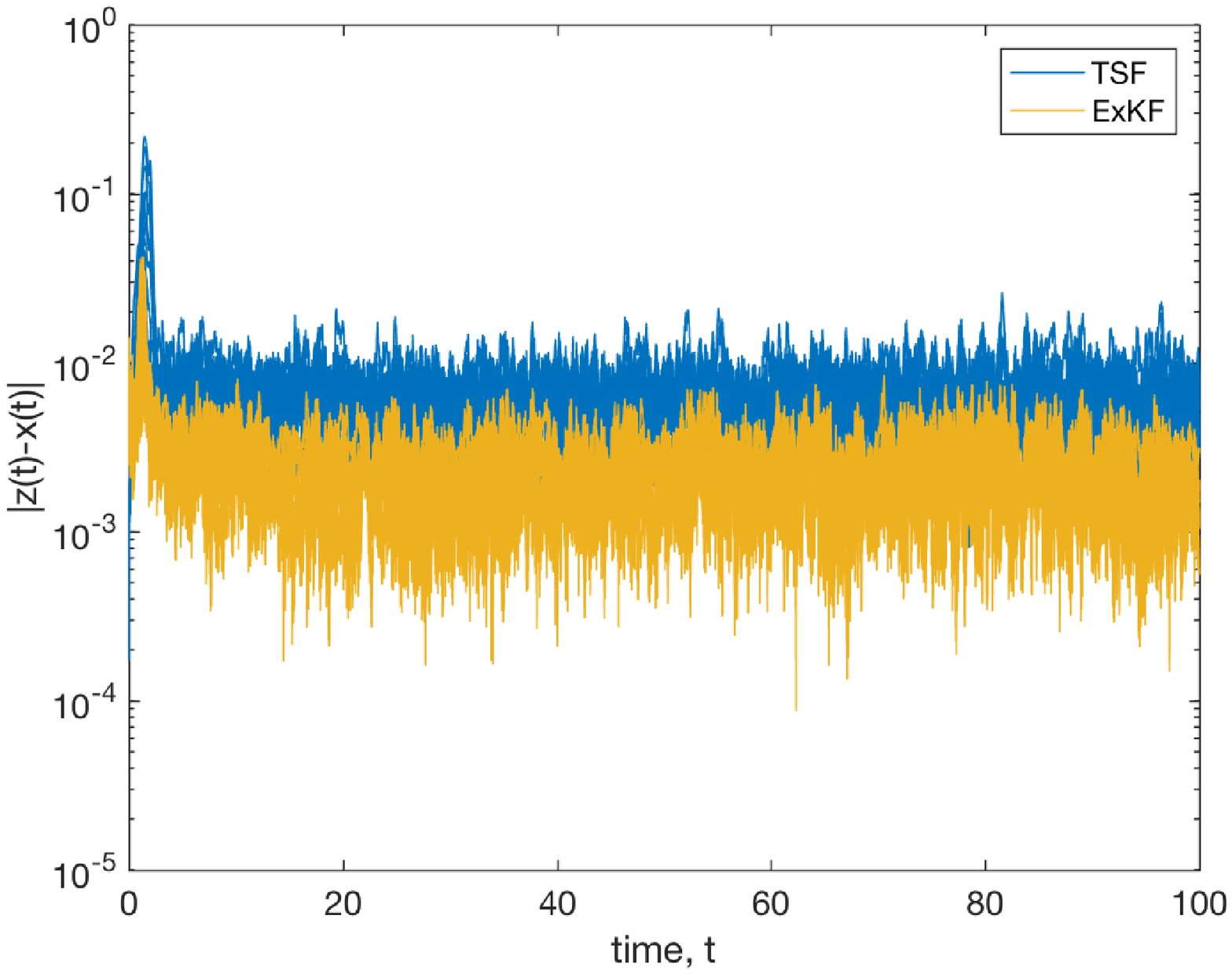}
\caption{Comparison of the filter \cref{eq:TSFilter} and the ExKF \cref{eq:linmnmx}
for the Lorenz '96 model \cref{eq:L96} with $k=8$.  Left, the errors $\|\xi(t)\|$ for a 10-member ensemble of perturbed initial conditions. Right, the errors $\|\xi(t) \|$  for a 10-member ensemble with random observational error.\label{fig7-8}}
\end{centering}
\end{figure}

\textbf{Comparisons with ExKF.} Recall from~\cite{KrenerMinimax} that the Kalman-Bucy filter for linear systems (with no model error), i.e. $f(t,m)=A(t)m$ is equivalent to the minimax filter:
\begin{equation}
  \label{eq:linmnmx}
  \begin{split}
      \dot m &= f(t,m) + PH^\intercal C(y-Hm), \quad m(0)=\est_0,\\
 \dot P &= \D{f}(t,m(t))P+P\D{f}(t,m(t))^\intercal(t) - PH^\intercal C H P\, \quad P(0) = P_0, \\
 \est_0^\intercal &P_0 \est_0 + \int_0^T \eta^\intercal(t) C(t) \eta(t) dt \le 1.
  \end{split}
\end{equation}
Instead of dealing with stochastic differential equations, for which the Kalman-Bucy filter is formulated, we stay within the (equivalent) deterministic framework. Namely, we take $\eta$ to be a measurable function satisfying the inequality in \cref{eq:linmnmx}.
We discretize~\cref{eq:linmnmx}  using the explicit fourth order Runge-Kutta method\footnote{A more appropriate integrator for \cref{eq:linmnmx} is an implicit symplectic integrator as detailed in \cite{ZhukCDC14,ZhukSISC13}. However for the sake of comparison we retain the explicit Runge-Kutta method here.}, reducing the
step size to $\Delta t=0.001$ to ensure stability on the time interval $T=100$. We begin with noise-free observations, and use an ensemble of 10 initial conditions, $\est_0=\truth(0)+0.01\eta$. In this case, to make sure that the error in the initial condition satisfies the inequality in \cref{eq:linmnmx}, we set $P(0)=\frac1{4 d\times 10^{-4}}I\approx 139 I$ and $C=I$, as $\|\sigma \eta\|^2$  is $\chi^2$-distributed with mean $k\cdot \sigma^2$. Clearly, large $P(0)$ and $C$ represent high trust in the initial condition and observations. As seen in the left panel of \cref{fig7-8}, the estimation error of the ExKF estimates is around $10^{-5}$ at the end of the interval, whereas the error of~\cref{eq:TSFilter} decays to machine precision.

\textbf{Noisy observations.} We next simulate a 10-member ensemble for which the initial condition is exact, $x_0 = z_0$, but the observations $y(t)$ are perturbed by random noise at each time step, $y(t) = H\truth(t)+0.01\eta$. Note that the expected value of the norm of the observational noise is given by:  $\mathbb{E} \, \sigma\|\eta\| \approx \frac{\sigma\sqrt{2 k-1}}{\sqrt{2}}$ for $k$ large. As demonstrated in the right panel of \cref{fig7-8},the estimation errors $\|\est-\truth\|$ of both the filter \cref{eq:TSFilter} and the ExKF \cref{eq:linmnmx} are less than the mean of the norm of the observational noise given by $\mathbb{E}\, \sigma\|\eta\| \approx 0.1275$, on the interval $(0,100)$. The ensemble mean estimation errors level off at around $10^{-2}$ for both methods, with the error of the \cref{eq:TSFilter} being approximately twice that of the ExKF \cref{eq:linmnmx}.
Recall from~\Cref{r:LQreconstruction} that this is as good as can be hoped for linear systems in the presence of noisy observations.

We stress that the observational noise introduces an additional non-linear term $p  Q \widetilde{Q}^{\intercal}  H^\intercal\eta$ in the error equation \cref{eq:error}; hence, large $p>0$ not only makes the discretized equation stiff but, importantly amplifies the noise! On the other hand, $\widetilde{Q}^{\intercal}$ acts as a projection onto the range of $H^\intercal HQ$, and thus $\widetilde{Q}^{\intercal}  H^\intercal\eta$ in fact represents the projection of $\eta$ onto the range of $HQ$. Hence, the norm of $Q \widetilde{Q}^{\intercal}  H^\intercal\eta $ is not increasing, and it may even be $0$ if $\eta$ is not in the range of $HQ$. The experiment shows that the amplification provided by $p=10$ is minor.


\subsection{Burgers equation}
\label{sec:burgers-equation}
As a second example, we discretize the Burgers(-Hopf) equation $u_t = - uu_x$ using the finite difference scheme:
\begin{equation}
  \label{eq:burgers} \dot{u}  = f(u), \qquad
\dot{u}_i = -\frac{1}{6 \Delta x} \left( u_i (u_{i+1} - u_{i-1}) + (u_{i+1}^2 - u_{i-1}^2), \right)
\end{equation}
again taken on a periodic lattice ($\Delta x = 2\pi/d$), which has the properties that (i) the quadratic energy $\sum_i u_i^2$ is conserved, implying that every sphere in $\Rs$ is invariant under the motion of the system and $\|u\|$ is constant, and (ii) $\operatorname{tr} (\D f(u)) = 0$, implying that the flow conserves the volume of the phase element.
Consequently, this equation is not dissipative in contrast to the $L96$ system, and so the effects of error in the initial condition, and the observational noise are expected to be more pronounced. On the other hand, it is unclear if this system possesses an invariant ergodic measure making it subject to the conditions of the Oseledec theorem, so our theory may not apply to this test case.

\begin{figure}[htpb]
\begin{centering}
\includegraphics[width=0.51\textwidth]{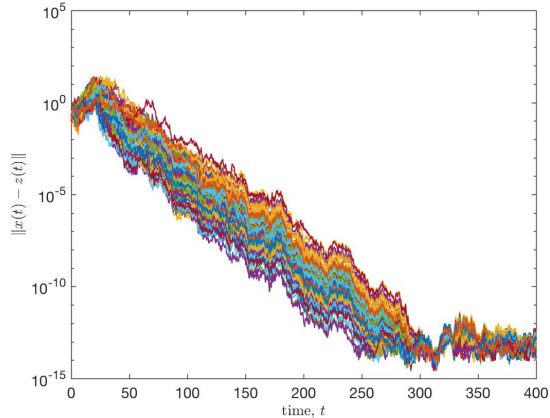}
\caption{Convergence of filter \cref{eq:TSFilter} for the discretized Burgers equation \cref{eq:burgers} with $k=11$, showing the errors $\|\xi(t)\|$ for a 100-member ensemble of perturbed initial conditions.\label{fig9}}
\end{centering}
\end{figure}

\textbf{Noise-free observations.} We set $d=18$ and take $\truth_i(0)\sim U(0,1)$, $i=1,\dots,d$.
For the filter~\cref{eq:TSFilter} we again take $H$ to be the eigenvectors of the discrete Laplace operator corresponding to the $k$ leading eigenvalues.
The filter~\cref{eq:TSFilter} was integrated by RK4 with $\Delta t=0.01$ and $p=20$ on the interval $t\in[0,400]$. For $k=9$ we found that $\dot X = \D{f}(\est(t))X$ possesses $10$ nonnegative exponents, and tangent dynamics $\dot \err = (A(t) - L(t,\est)H^\intercal(t) H(t))\err$ has $1$ nonnegative exponent. Taking $k=10$ we observed convergence, i.e.~$\|\est-\truth\|\le 10^{-14}$, provided $\est_0=\truth(0)+0.0005\eta$.
This suggests that the basin of attraction of the trivial solution of~\cref{eq:error} is rather small. However, for $k=11$ this basin increases significantly: the exponential decay of the estimation error for an ensemble of 100 initial conditions $\est_0=\truth(0)+0.01\eta$
(relative error $\frac{\|\truth_0-\est_0\|}{\|\truth_0\|}$ is up to $45\%$) is evident in \cref{fig9}.

\begin{figure}[h]
\begin{centering}
\includegraphics[width=0.49\textwidth]{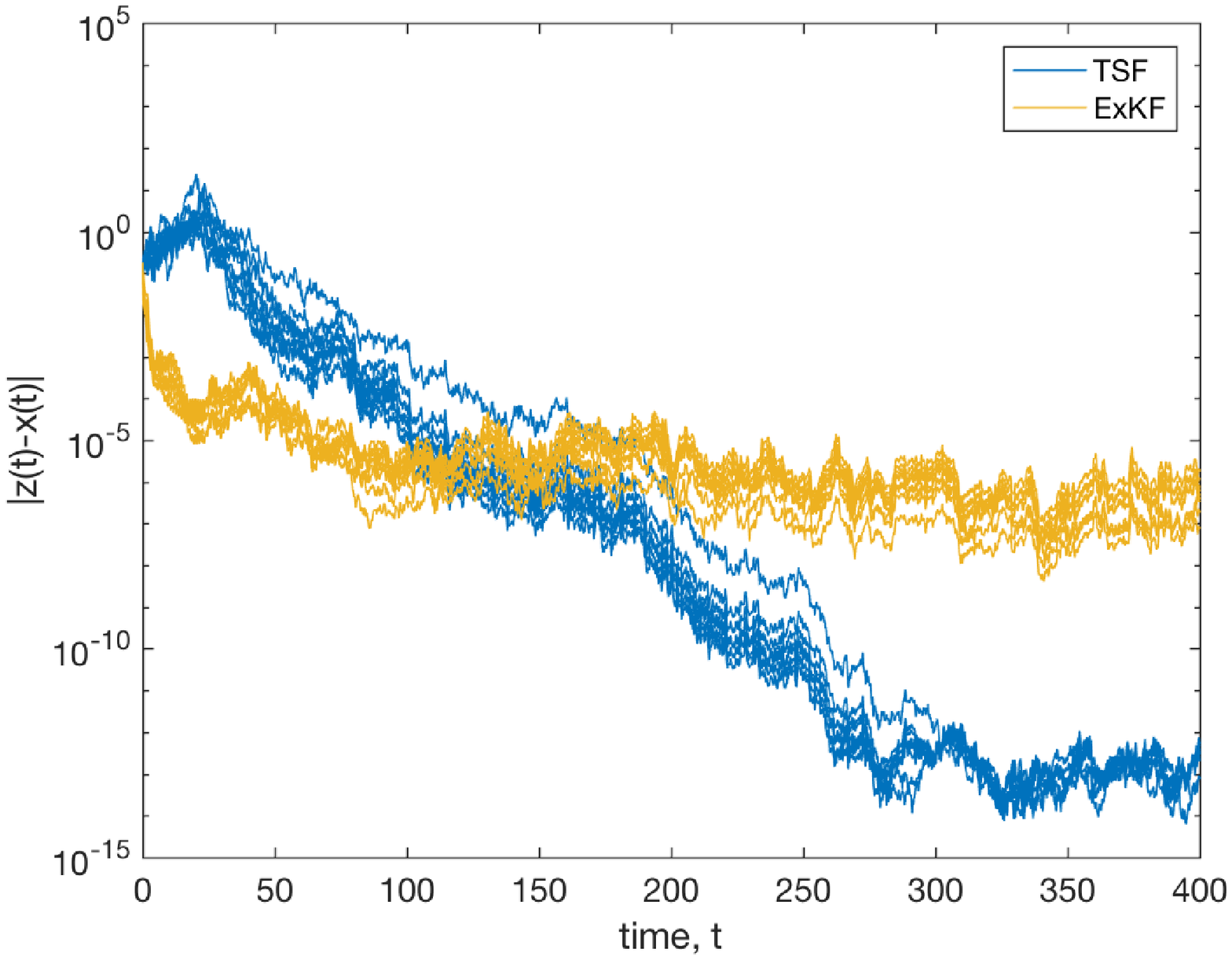}
\includegraphics[width=0.49\textwidth]{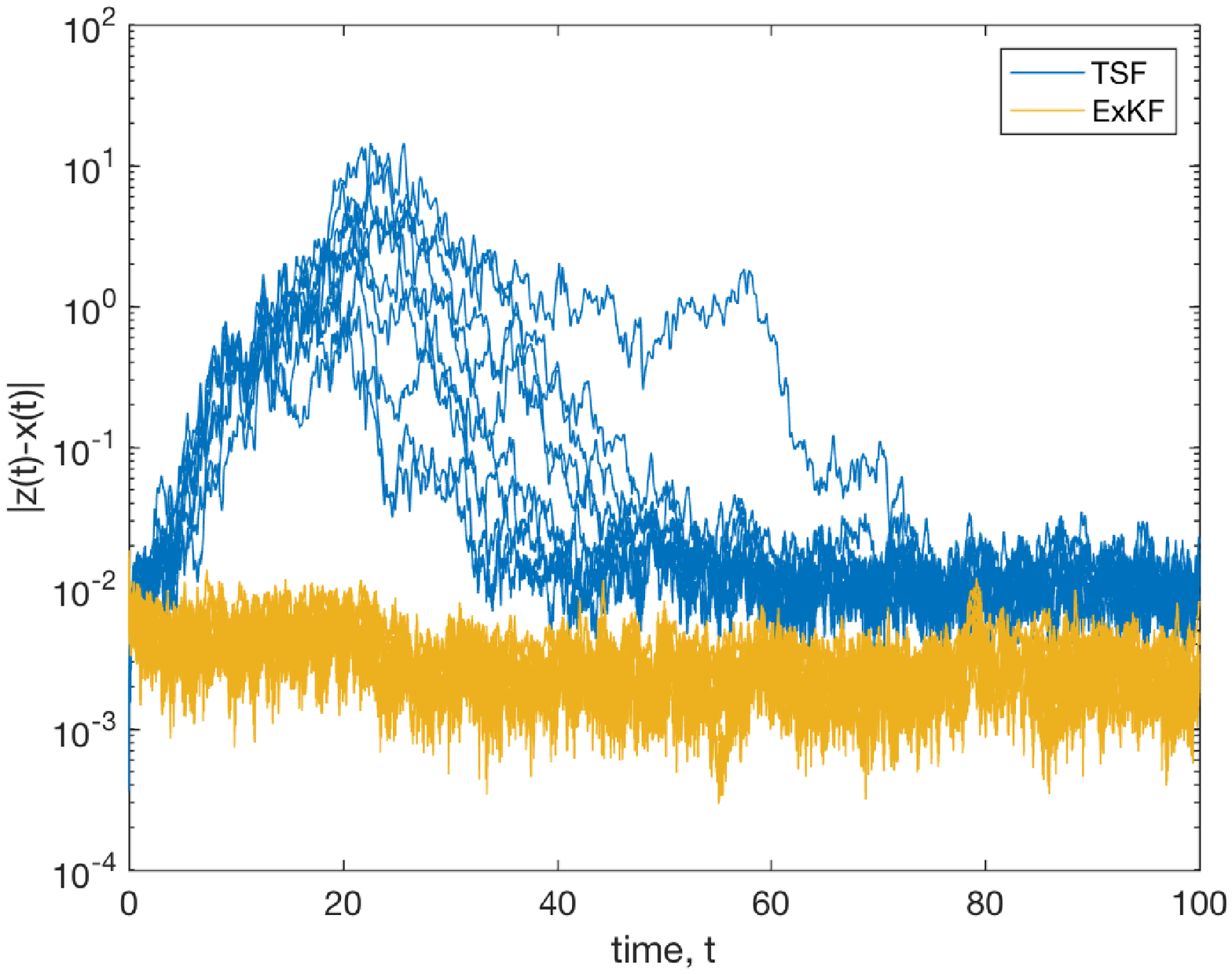}
\caption{Comparison of the filter \cref{eq:TSFilter} and the ExKF \cref{eq:linmnmx}
for the discretized Burgers equation \cref{eq:burgers} with $k=11$.  Left, the errors $\|\xi(t)\|$ for a 10-member ensemble of perturbed initial conditions.  Right, the errors $\|\xi(t) \|$  for a 10-member ensemble with random observational error.\label{fig10-11}}
\end{centering}
\end{figure}

\textbf{Comparisons with ExKF.}
To compare with the ExKF method \cref{eq:linmnmx} we again reduce the
step size to $\Delta t=0.001$. We begin with noise-free observations, and use an ensemble of 10 initial conditions, $\est_0=\truth(0)+0.01\eta$. In this case, to make sure that the error in the initial condition satisfies the inequality in \cref{eq:linmnmx}, we set $P(0)=\frac1{4 d\times 10^{-4}}I\approx 139 I$ and $C=I$, as $\|\sigma \eta\|^2$  is $\chi^2$-distributed with mean $k\cdot \sigma^2$. As seen in the left panel of \cref{fig10-11}, the error of the ExKF estimates is again around $10^{-5}$ at the end of the interval whereas the error of~\cref{eq:TSFilter} decays to machine precision.

\textbf{Noisy observations.} We simulate a 10-member ensemble for which the initial condition is exact, $x_0 = z_0$, but the observations $y(t)$ are perturbed by random noise at each time step, $y(t) = H\truth(t)+0.01\eta$.  As shown in the right panel of \cref{fig10-11} the ensemble mean estimation error levels off at around $2\times 10^{-2}$. Here the mean of the norm of the observational noise is $0.033$.  The convergence of the filter \cref{eq:TSFilter} is irregular at the beginning of the interval, probably due to slow convergence of $Q(t)$ to the basis of the nonstable tangent space from its random initial condition.  The error in the ExKF method is ultimately smaller than that of the filter \cref{eq:TSFilter} by a factor 3.



\section{Concluding remarks}
\label{sec:concluding-remarks}
In this paper we have drawn an explicit connection between the dimension of the nonstable tangent space of a continuous dynamical system, as quantified by the number of nonnegative Lyapunov exponents, and the necessary dimension of the (time-variant) observation operator $H$ of a sequential data assimilation process. We formulated a detectability condition that, when satisfied, provides a necessary and sufficient condition for convergence of the new filter~\cref{eq:TSFilter} that utilizes an explicit partition of the tangent space into stable and nonstable subspaces. The new filter is comparable to the Extended Kalman Filter for perturbed initial conditions and noisy observations, and appears to be more robust in the sense that convergence is observed for an order of magnitude larger than the time step.


\appendix

\section{Proofs}
\label{sec:proofs}

\begin{proof}[Proof of~\cref{l:luenberger}]

The first statement follows directly from~\cref{d:detect_LTI} and~\cref{LEfun}. Let us prove the second statement. Recall that $\Rs=\ker O^s\oplus \ker^\perp O^s$ and take $w\in \ker^\perp O^s$. \Cref{d:detect_LTI} implies that $\int_0^T \| H e^{At}w\|dt>0$ for any $T>0$. Hence, $H\fatQ_i(t)\ne 0$ for any $\fatQ_i$ from the non-stable tangent subspace of $\dot X = A X$; see~\cref{eq:LEQA}. The latter shows that \Cref{d:detect_LTI} implies \Cref{d:detect-At}, which is equivalent to the existence of the observer by~\Cref{p:gain}.

Let us prove the last statement. Take $w\in \ker O^s$ and assume that $e^{At}w$ does not decay, and, for some $L$, the observer exists. Note that $y(t) = H  e^{At}w = 0$. Hence, by~\Cref{d:luenberger}, the solution of~\cref{eq:Luenberger} converges to $0$ for any $\est_0$ as in this case~\cref{eq:Luenberger} coincides with the error equation $\dot\err=(A-LH)\err$, and, by~\Cref{d:luenberger}, the solution of the latter decays to $0$ exponentially fast for any $\err(0)$. But this contradicts the orinal assumption ($e^{At}w$ does not decay!). This contradiction proves the last statement of the lemma.

Let us also prove that $\ker(W(0,T)) = \ker(O^s)$. Indeed, \[
W(0,T)v=0\Leftrightarrow \int_0^T \|He^{At}v\|^2_{R^p}dt = 0 \Leftrightarrow \|He^{At}v\|^2_{R^p}\equiv 0\,,
\]
and since $t\mapsto He^{At}v$ is a smooth function, it follows that $He^{At}v=0$ for any $t\in[0,T]$. By differentiating the latter equality we find that $HA^je^{At}v=0$ for any $j\ge0$ and $t\in[0,T]$. Setting $t=0$ we get that $W(0,T)v=0$ implies $O^sv=0$. On the contrary, by definition of $s$ we have that \[
(A^{s+p})^\intercal H^\intercal\in \operatorname{span}\{H^\intercal\,, A^\intercal H^\intercal,\dots,(A^\intercal)^s H^\intercal\}\,,\quad p\ge 1\,,
\]
Hence, $O^sv=0$ implies that $v^\intercal (O^s)^\intercal = 0$ and so $v^\intercal (A^j)^\intercal H^\intercal=0$ for any $j\ge0$. Therefore, $v^\intercal e^{A^\intercal t}H^\intercal=0$.
\end{proof}

\begin{proof}[Proof of~\cref{p:gain}]
Note that $\dot \err = (A(t) - L(t) H(t)) \err$, and $A$, $H$ are bounded matrix-valued functions by~\cref{eq:AH}. Hence, $\err$ decays to $0$ exponentially fast iff  the LEs $\mu_1 \ge \cdots \ge \mu_d$ of
\begin{equation}\label{eq:W}
\dot W = (A(t) - L(t) H(t))W\,,
\end{equation}
satisfy $\mu_1 < -\kappa<0$ (see~\cite[p.6]{BaPe02}). We stress that the forward regularity is not required for the latter statement to hold as per~\cref{eq:exp_decay}. However, it becomes important when we invoke~\cref{d:detect-At} to prove that $\err$ decays to $0$ exponentially fast. Our proof is based on the following simple observation: if $X$ solves~\cref{fundmat}, $L$ satisfies~\cref{eq:gainL}, and~\cref{eq:forward_reg} holds, then:
\begin{enumerate}[label= (U), align=left, leftmargin=*]
\item \label{uu}
if $X(t)=\fatQ(t)\fatR(t)$ and $\fatR_1 = \fatQ^\intercal W$ then $\fatR_1(t)$ is upper-triangular with positive diagonal, i.e.~$\fatQ\fatR_1$ represents the unique QR-decomposition of $W$.
\end{enumerate}
Assume for now that \ref{uu} holds true. Recall from~\cref{eq:LE_QAQ} that, given the unique QR-decomposi\-tion of $X$, i.e.~$X=\fatQ\fatR$, one can compute the $i$th Lyapunov exponent of~\cref{fundmat}, $\lambda_i$ by evaluating the limit of the quantity $\frac 1t \int_0^t Q_i^\intercal A Q_ids$, which depends only on $A$ and the $i$th column of $\fatQ$, $Q_i$. Hence, by \ref{uu}, the $i$th Lyapunov exponent of~\cref{eq:W}, $\mu_i$, depends only on $A-LH$ and the same $Q_i$: \[
\mu_i = \lim_{t\to+\infty} \frac 1t \int_0^t Q_i^\intercal (A - LH) Q_ids = \lambda_i(X) - \lim_{t\to+\infty} \frac 1t \int_0^t Q_i^\intercal LH Q_ids\,.
\]
Now, by \cref{eq:gainL} it follows that $LH=pQ\widetilde{Q}^\intercal H^\intercal H$, and hence, by~\cref{eq:fat-skin-Q}, we get:
\begin{equation}
  \label{eq:QLHQ}
  \fatQ^\intercal L H \fatQ  = \begin{bmatrix} Q^\intercal  \\ Q_\perp^\intercal \end{bmatrix}
	p\, Q\, \widetilde{Q}^\intercal H^\intercal H \begin{bmatrix} Q & Q_\perp \end{bmatrix}
	= \begin{bmatrix}  p\, \widetilde R  & p\, \widetilde{Q}^\intercal H^\intercal H Q_\perp \\ 0&	0\end{bmatrix}
\end{equation}
Clearly, $Q_i^\intercal LH Q_i = p\widetilde R_{ii}$, $i=1,\dots,k$, and so: \[
\mu_i = \lambda_i - p \lim_{t\to+\infty} \frac 1t \int_0^t \widetilde R_{ii} ds\,.
\]
Since $\widetilde R_{ii}\ge0$, it follows by the previous equality that $\mu_i<0$ for $i>k$, and that the $k$ leading LEs $\mu_1 \ge \cdots \ge \mu_k$ of
\cref{eq:W} are negative iff \[
\exists p>0:\quad p \lim_{t\to+\infty} \frac 1t \int_0^t \widetilde R_{ii} ds\,>\lambda_i>0\,,\quad i=1,\dots,k
\]
But this is the case if and only if $\matpair{A}{H}$ is detectable in the sense of~\cref{d:detect-At}. Note that selecting \[
	p> \frac{\kappa+\max_j\{\lambda_j\}}{\min_{j}\{\lim_{t\to +\infty} \frac1t\int_0^t \widetilde R_{jj}(t) dt\}}
\] we achieve the desired inequality, namely $\mu_1<\kappa$. This completes the proposition's proof.

Let us now prove \ref{uu}. Let $X$ solve~\cref{fundmat}, $X(0)=\fatQ_0\,\fatR_0$, and $\fatQ\in\R^{d\times d}$, $\fatR\in\R^{d\times d}$ denote the unique solution of~\cref{eq:QR_R,eq:QR_Q}, and $Q\in\R^{d\times k}$, $R\in\R^{k\times k}$ solve~\cref{qr1,qr2} for $k\le d$. Define $\fatR_1(t):= \fatQ(t)^{\intercal} W(t)$ and assume that $W(0)=X(0)$. Then
\[
	\dot{\fatR}_1 = \dot{\fatQ}^{\intercal} W + \fatQ^{\intercal} \dot{W} =
	(\fatQ^\intercal A \fatQ - \fatS - \fatQ^\intercal L(t) H \fatQ) \fatR_1\,, \quad \fatR_1(0)=\fatR_0\,.
\]
Here $\fatS$ denotes $S$ defined in~\cref{eq:QR_Q}. Recall from~\cref{eq:QR_R} that $\fatB = \fatQ^\intercal A \fatQ - \fatS$ is upper-triangular.  Furthermore, by~\cref{eq:QLHQ}, $\fatQ^\intercal L(t) H \fatQ$
is also upper-triangular. Consequently, $\fatR_1(t)$ is
upper-triangular with positive main diagonal if it is so initially. Thus $\fatR_1$ and $\fatQ$ verify~\ref{QR2}, and hence, by~\ref{QR1}, $\fatR_1(t)\fatQ(t)$ coincides with the unique QR-decomposition of $W$.
\end{proof}
\begin{proof}[Proof of~\cref{t:tsf}]
We first prove \ref{s2} $\Rightarrow$ \ref{s1}. Note that $\dot \err = \dot\truth-\dot\est = f(t,\truth)-f(t,\est)-L(t,\est)H(t)\err(t)$. Recall that $f = (f_1,\dots,f_d)^\intercal$, and by \ref{a1} and \ref{a2} of \cref{a:lindyn} $f$ has bounded Jacobian (w.r.t.~$x$) $A(t)=\D{f}(t,\est(t))$ and $N_i=D^2 f_i(t,x)$, the Hessian of each $f_i$ w.r.t.~$x$ is bounded on every compact set $K$, uniformly w.r.t.~time $t>0$. Moreover, by \ref{a2} we have that:
\begin{equation}
  \label{genLipschitz}
\|N(t,\xi_1,\est(t))-N(t,\xi_2,\est(t))\|\le C(K)\|\xi_1-\xi_2\|^2\,, \quad C(K):=\max_i C_i(K)
\end{equation}
is uniformly H\"{o}lder continuous with exponent 2 on every compact set $K\subset\Rs$. Define $g(s):=f_i(s \err +\est)$. Then $g(1)-g(0)=\int_0^1 g'(s) ds$, or, equivalently \[
f_i(t,\err(t) +\est(t)) - f_i(t,\est(t)) = \int_0^1 \err^\intercal(t) \nabla f_i (t,s \err(t) + \est(t)) ds
\] Now, define $g(\tau;s,t):=\err^\intercal(t)\nabla f_i (t,\tau s \err(t) + \est(t))$. We get that \[
  f_i(t,\err(t) +\est(t)) - f_i(t,\est(t)) - \err^\intercal(t)\nabla f_i (t,\est(t))
=  \int_0^1 (g(1;s,t) - g(0;s,t)) ds
= N_i(t,\err(t),\est(t))
\]
This implies that \[
 f(t,\truth)-f(t,\est) = f(t,\err+\est)-f(t,\est) = \D{f}(t,\est) + N(t,\err,\est)
\] and so the estimation error $\err$ solves~\cref{eq:error}. Now, \ref{s1} follows\footnote{See (1.4.14) on p.~29 of~\cite{BaPe02}} from~\cite[p.27, Thm.~1.4.1]{BaPe02} and~\cite[p.~29, Thm.~1.4.3]{BaPe02}.

We next prove \ref{s1} $\Rightarrow$ \ref{s2}. Note that~\cref{eq:error} is equivalent to the following integral equation: \[
\err(t) = X(t) X^{-1}(0)\err_0 + \int_0^t X(t)X^{-1}(s) N(s,\err(s),\est(s)) ds
\] provided $\dot X = (A-LH)X$, $X(0)=X_0$. Take some $q$ and consider a linear equation $\dot v = (A-LH)v + q$, $v(0)=\err_0$ or its integral representation: \[
v(t) =  X(t) X^{-1}(0)e_0 + \int_0^t X(t)X^{-1}(s) q(s) ds
\] Since $X$ is independent of $\err$, it follows that $v(t)=\err(t)$ provided $q(s) = N(s,\err(s,\err_0),\est(s))$. Hence, \ref{s2} is verified if we can show that $X(t) X^{-1}(0)\err_0$ decays to $0$ exponentially fast for any $\err_0$: $\|\err_0\|<\varepsilon$, i.e.\ all the LEs $\mu_d\le\dots\le\mu_1$ of $\dot X = (A-LH)X$ are negative. Let $X=\fatQ\fatR$ be the (unique) full QR decomposition of $X\in\R^{d\times d}$, and set $z:=\fatQ^\intercal v$, $p:=\fatQ^\intercal q$. Then $\dot z = B(t) z + p$, $z(0) = \fatQ^\intercal(0)\err_0$ and $\|z\|=\|v\|$. Here $B$ is an upper-triangular matrix defined as in~\cref{eq:QR_R}, and $\mu_i=\lim_{t\to+\infty} \frac 1t \int_0^t B_{ii}(s) ds$ by~\cref{eq:LE_QAQ}. Note that $\dot z_d = B_{dd} z_d + p_d$, or, equivalently: \[
z_d(t) = e^{\int_0^t B_{dd}(s) ds} z_d(0) + e^{\int_0^t B_{dd}(s) ds}\int_0^t e^{-\int_0^s B_{dd}(\tau) d\tau} p_d(s) ds\,.
\]
Assume that $\mu_d>0$. Then $\frac 1t\int_0^t B_{dd}(s) ds> \mu_d-\delta$ for a small $\delta>0$ and all $t>t^\star$. Thus \[
-\int_0^s B_{dd}(\tau) d\tau \le 0 \Rightarrow \gamma:=\int_0^{+\infty} e^{-\int_0^s B_{dd}(\tau) d\tau} p_d(s) ds <+\infty
\] as $p(s) = \fatQ^\intercal (s) N(s,\err(s,\err_0),\est(s))$, and $N$ satisfies~\cref{genLipschitz}, and $\|\err(t)\|\le C e^{-at}\|\err_0\|$ by \ref{s1} Let \[
\theta(t):=\gamma - \int_0^t e^{-\int_0^s B_{dd}(\tau) d\tau} p_d(s) ds\,.
\] Then $\|\theta(t)\|\le\varepsilon$ provided $t>t^\star_1$. Without loss of generality we can assume that $z_d(0)<\gamma$. But then, for $t>t^\star_1$ we have that \[
z_d(t)=e^{\int_0^t B_{dd}(s) ds} \left(z_d(0) + \gamma - \theta(t)\right)\to +\infty\,,
\] as $\int_0^t B_{dd}(s) ds> t(\mu_d-\delta)$ provided $t>t^\star$ and $\mu_d-\delta>0$ so that $e^{\int_0^t B_{dd}(s) ds}$ grows unbounded. This condradicts $2.$ as $|z_d|\le \|z\|=\|x\|=\|\err\|\le C e^{-at}\|\err_0\|$. Now, if $\mu_d=0$ then, for any $\delta$ such that $0<\delta<a$ there exists $t_\delta>0$ such that $ -\delta t< \int_0^t B_{ii}(s) ds<\delta t$, for all $t>t_\delta$. Hence, $\gamma$ is still bounded as $|p_d|<\|N\|\le \tilde C\|\err_0\| e^{-at}$. Hence, $z_d>e^{-a t}\tilde C_1$ which again contradicts $2.$ Hence $\mu_d<0$. Noting that $\dot z_{d-1} = B_{d-1d-1}z_{d-1} + B_{d-1 d} z_d + p_{d-1}$, and that $B_{d-1 d}$ is proportional to the function $s\mapsto Q_{d-1}^\intercal(s) A(s) Q_d(s)$ which is bounded from above as $\|A\|<+\infty$ we can rewrite the equation for $z_{d-1}$ as follows: \[
z_{d-1}(t) = e^{\int_0^t B_{d-1 d-1}(s) ds} z_{d-1}(0) + e^{\int_0^t B_{d-1d-1}(s) ds}\int_0^t e^{-\int_0^s B_{d-1d-1}(\tau) d\tau} \tilde p_{d-1}(s) ds\,,
\]
where $\tilde p_{d-1} = p_{d-1}+B_{d-1 d} z_d$ decays to zero exponentially fast. Then it is easy to demonstrate that $\mu_{d-1}=\lim_{t\to+\infty} \frac 1t \int_0^t B_{d-1d-1}(s) ds<0$ by the same argument as was used above to show that $\mu_d<0$. Repeating this argument for every $d-i$, $i<d$ we obtain that indeed \ref{s1} $\Rightarrow$ \ref{s2}. This completes the proof.
\end{proof}
\begin{proof}[Proof of~\cref{c:tsf}]
In the proof of~\cref{p:gain} we demonstrated that the equation $\dot e = (A(t) - L(t) H(t)) e$ has negative LEs iff $\matpair{A}{H}$ is detectable, and the gain $L$ is defined by~\cref{eq:gainL}. Note that the aforementioned statement about LEs is exactly the statement \ref{s2} of~\cref{t:tsf} which is equivalent to \ref{s1}. This completes the proof.
\end{proof}

\bibliographystyle{abbrv}

\end{document}